\documentclass[12pt, letterpaper]{article}
\usepackage{geometry}
\usepackage{graphicx}
\usepackage{amssymb}
\usepackage{epstopdf}
\usepackage{amsmath}
\usepackage{amsmath, ams cd}
\usepackage{slashbox}
\pdfoutput=1

\newtheorem{prop}{Proposition}[section]
\newtheorem{theorem}[prop]{Theorem}
\newtheorem{lemma}[prop]{Lemma}

\newtheorem{remark}[prop]{Remark}

\begin{document}

\title{Virtually fibred Montesinos links of type $\widetilde{SL_2}$}

\author{Xiao Guo\thanks{E-mail: xiaoguo@buffalo.edu}, Yu Zhang\thanks{E-mail: yz26@buffalo.edu}\\Department of Mathematics, SUNY at Buffalo}

\date{ } 

\maketitle

\begin{abstract}
We find a larger class of virtually fibred classic Montesinos links
of type $\widetilde{SL_2}$, extending a result of Agol, Boyer and
Zhang.
\end{abstract}

\section{Introduction}

A $3$-manifold is called \textit{virtually fibred} if it has a finite cover
which is a surface bundle over the circle. A link in a connected
$3$-manifold is said to be \textit{virtually fibred} if its exterior is a
virtually fibred $3$-manifold. Thurston conjectured that all closed
hyperbolic $3$-manifolds and all hyperbolic links in closed
$3$-manifolds are virtually fibred. This conjecture, which has been
named as \textit{virtually fibred conjecture}, is one of the most fundamental
and difficult problems in $3$-manifold topology.

\begin{figure}
\begin{center}
\includegraphics{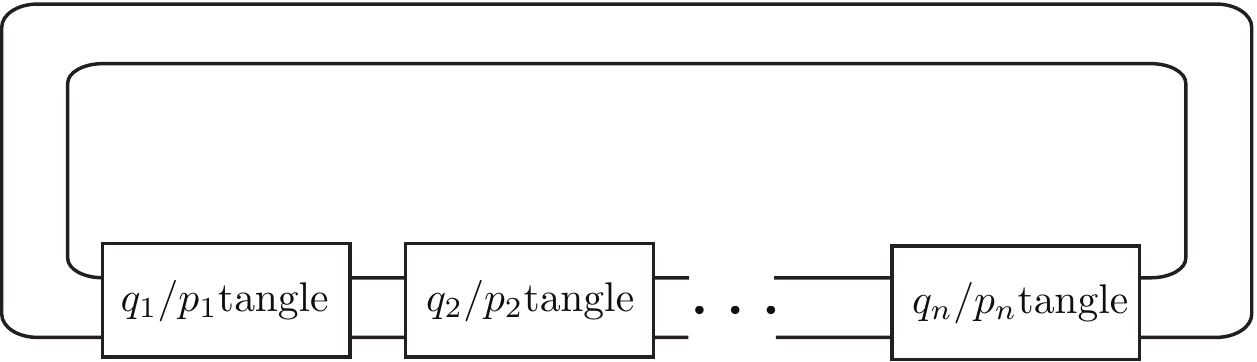}
\end{center}
\caption{\label{cml} Classic Montesinos link}
\end{figure}

Recall that a link $K$ in $S^3$ is called a {\it generalized
Montesinos link} if the double branched cover $W_K$ of $(S^3, K)$ is
a Seifert fibred $3$-manifold. Such a link $K$ is further said to be
\textit{of type $\widetilde{SL_2}$} if the canonical geometric structure on
$W_K$ is from the $\widetilde{SL_2}$-geometry. When every component
of the branched set in $W_K$ is not a fiber of the Seifert fibration
of $W_K$, $K$ is called a \textit{classic Montesinos link}. Recent
work of Walsh [Wa] and work of Agol-Boyer-Zhang [ABZ] combined
together solved the virtually fibred conjecture for all generalized
Montesinos links in $S^3$ which are not classic Montesinos links of
type $\widetilde{SL_2}$. Agol-Boyer-Zhang [ABZ] also gave an
infinite family of virtually fibred classic Montesinos links of type
$\widetilde{SL_2}$. In this paper, we extend the latter result of
[ABZ] to a larger family of classic Montesinos links of type
$\widetilde{SL_2}$. Note that every classic Montesinos link has a
cyclic rational tangle decomposition of the form
$(q_1/p_1,q_2/p_2,...,q_n/p_n)$ with all  $p_i\geq 2$ as shown in
Figure \ref{cml}. We prove

\begin{theorem}\label{main}
If a classic Montesinos link $K$ has a cyclic rational tangle
decomposition of the form $(\displaystyle{\frac{q_1}{p},
\frac{q_2}{p}, \cdots, \frac{q_n}{p}})$ with  $p\geq 3$ odd,
then $K$ is virtually fibred.
\end{theorem}

This theorem is proved in [ABZ] when $n$ is a multiple of $p$. Our
approach follows closely to that of [ABZ]. At a number of places we
need to deal with some new issues that arise. We shall describe our
ways of dealing with these issues when we reach these places.

Let $K$ be a link as given in Theorem 1.1. The base orbifold ${\cal
B}_K$ of the Seifert fibred space $W_K$ is a $2$-sphere with $n$
cone points each has order $p$.  Let $f:W_K\rightarrow \mathcal{B}_K$
be the Seifert fibration which is invariant under the covering
involution $\tau:W_K\rightarrow W_K$. The Euler number $e(W_K)$ of
$W_K$ and the orbifold Euler characteristic $\chi({\cal B}_K)$ of
${\cal B}_K$ are given by the following formulas:
$$e(W_K)=-\frac{\sum_{i=1}^nq_i}{p},\;\;\; \chi({\cal B}_K)=2-n+
\frac{n}{p}.$$ Note that $W_K$ has the $\widetilde{SL_2}$-geometry
precisely when $e(W_K)\ne 0$ and $\chi({\cal B}_K)<0$. Thus $K$ is
of type $\widetilde{SL_2}$ precisely when

{\bf Case (1)}. $n=3$, $p\geq 5$, and $q_1+q_2+q_3\ne 0$, or

{\bf Case (2)}. $n>3$ and $q_1+...+q_n\ne 0$.

We shall split our proof of Theorem \ref{main} into these two cases, given
in  Section 2 and Section 3 respectively. All main ingredients of the
proof will occur already in Case (1); the proof of Case (2) will be
a quick  generalization.

In the remainder of this section, we shall give a few additional
notes and notations which will be used throughout this paper.
 From now
on we assume that $K=K(q_1/p,...,q_n/p)$ is a link satisfying the
conditions of Case (1) or (2) listed above.

First note that $K$ is a single component knot if $q_1+...+q_n$ is
odd, and a two component link otherwise. The double branched cover
$W_K$ has the $\widetilde{SL_2}$-geometry and its base orbifold
${\cal B}_K$  is hyperbolic. Let $\widetilde{K}$ be the corresponding
branched set in $W_K$ and let $K^*=f(\widetilde{K})$, where $f:W_K\rightarrow {\cal B}_K$ is the Seifert quotient map. By the orbifold theorem we
may assume that the covering involution $\tau$ on $W_K$ is an
isometry and the branch set $\widetilde{K}=Fix(\tau)$ is a geodesic
which is orthogonal to the Seifert  fibres by Lemma 2.1 of [ABZ].
Also $K^*$ is a geodesic in ${\cal B}_K$, and is an equator of
${\cal B}_K$ containing all the $n$ cone points which we denote by
$c_1, c_2, \cdots, c_n$ (indexed so that they appear consecutively
along $K^*$). The order of each cone point $c_i$ is $p$. Since $p$
is  odd, the  restriction map  of $f$ to each component of
$\widetilde{K}$ is a double covering onto $K^*$.

\section{Proof of Theorem \ref{main} in Case (1).}

In this section we prove Theorem \ref{main} in the case $n=3$, $p\geq 5$
and $q_1+q_2+q_3\ne 0$. Throughout this section, without further notice, $i, j, k$ as well as numbers expressed by them will be considered as non-negative integers mod $p$,
and $i\neq k$ will also be assumed.

\subsection{$\textbf{$K$}$ is a knot.}

Here is an  outline of the proof.

\noindent{\bf Step 1}\ \ Take  a specific $p^2$-fold orbifold cover $F$
of ${\cal B}_K=S^2(p,p,p)$ such that $F$ is a smooth surface. There
is a corresponding free cover $\Psi : Y \rightarrow W_K$ of the same
degree such that the base orbifold of the Seifert space $Y$ is $F$
and such that the following diagram commutes
\begin{equation*}
\begin{CD}
Y       @>{\hat{f}}>>      F\\
@VV{\Psi}V        @VV{\psi}V\\
W_K      @>{f}>>{\mathcal{B}_K}
\end{CD}
\end{equation*}
where $\hat f$ is the Seifert quotient map. Note that  $Y$ is a
locally-trivial circle bundle over $F$ since $F$ is a smooth
surface.

 Let $L= \Psi^{-1}(\widetilde{K})$. Then $L$ has exactly
$p^2$ components, which we denote by  $\{L_{i,j} : 1\leqslant i,
j\leqslant p\}$. Let $L_{i,j}^*=\hat{f}(L_{i,j})$. Then $\{L_{j,j}^{*} :
1\leqslant j\leqslant p\}$ are $p$ mutually disjoint simple closed
geodesics in $F$.

\noindent\textbf{Step 2}\ \ Construct a surface semi-bundle structure on
$M=Y - \overset{p}{\underset{j=1}{\cup}}\overset{\circ}{N}(L_{j,j})$,
where $N(L_{j,j})$ is a small regular neighborhood of $L_{j,j}$ in
$Y$, $1\leqslant j\leqslant p$.

Note that $M$ is a graph manifold with non-empty boundary, with
vertices $M_2^j $ and $Y_1$, where $M_2^{j}=
\hat{f}^{-1}(L_{j,j}^{*})\times[-\epsilon, \epsilon]-\overset{\circ}{N}(L_{j,j})$,
$1\leqslant j\leqslant p$.
$Y_1=M-\overset{p}{\underset{j=1}{\cup}} \overset{\circ}{M^j_2}$. We
construct a surface semi-bundle structure on $M$ by using
\cite{wy} and following \cite{abz}.

\noindent\textbf{Step 3}\ \ Isotope all $L_{i,k}$, $i\ne k$,
 such that they are transverse to the surface bundle in $M_2^j$,
 denoted by $\mathcal{F}_2^j$.
This is one of the key parts of the proof. Let
$U_{i,k}=\hat{f}^{-1}(L_{i,k}^{*})$, the vertical torus over
$L_{i,k}^{*}$. By the construction in Step 2, there are exactly two
singular points in the induced foliation by ${\cal F}_2^j$ on  every
component of $U_{i,k}\cap M_2^j$. 
We show that $L$ can be oriented such that the arcs $L_{i,k}\cap M_2^j$ whose images intersect in $F$, all travel from one of $\hat{f}^{-1}(L_{j,j})^{*}\times\{-\epsilon\}$ and $\hat{f}^{-1}(L_{j,j})^{*}\times\{\epsilon\}$ to the other. This property allows us to arrange these arcs in
$U_{i,k}\cap M_2^j$ such that  they always travel from the ``$-$"
side of the surface bundle in $M_2^j$ to the ``$+$" side.

\noindent\textbf{Step 4}\ \ Construct a double cover of $M$, denoted by
$\breve{M}$, so that $\breve{M}$ has a surface bundle structure. Denote the
corresponding double cover of $Y$ as $\breve{Y}$,the lift of $Y_1$ as $\breve{Y}_{1,1}\cup \breve{Y}_{1,2}$, and the lift
of $L_{i,k}$ as $\breve{L}_{i,k,1}\cup\breve{L}_{i,k,2}$.
 Perform certain Dehn twist operations on the
surface bundles of $\breve{Y}_{1,1}\cup \breve{Y}_{1,2}$ so that the $\breve{L}_{i,k,s}$ are transverse to the new surface bundles, $s=1,2$.

This is another key part of the proof. We perform these Dehn twist
operations along a union of vertical tori, $\Gamma \subset
\breve{Y}_{1,1}\cup \breve{Y}_{1,2}$. $\Gamma$ contains some
boundary parallel tori of $\breve{Y}_{1,1} \cup \breve{Y}_{1,2}$ as
in \cite{abz}, different form \cite{abz}, these boundary parallel tori intersect some arcs
of $\breve{L}_{i,k,s}\cap \breve{Y}_{1,s}$ two times with different direction or do not intersect at all, $s=1,2$. We call such arcs of
$\breve{L}_{i,k,s}\cap  \breve{Y}_{1,s}$ ``bad'' arcs. We
construct four specific  extra tori as members of $\Gamma$ to deal
with the ``bad'' arcs, $s=1,2$.

By Steps 1-4, the exterior of the inverse image of $\widetilde{K}$ in
$\breve Y$ has a surface bundle structure. It is a free cover of
the exterior of $\widetilde{K}$ in $W_K$, which in turn is a free double
cover of the exterior of $K$ in $S^3$. Thus $K$ is virtually fibred in $S^3$.

Now we fill in the details.

\noindent\textbf{Step 1}

Instead of taking the $p$-fold cyclic cover of $\mathcal{B}_K$ as in \cite{abz}, we construct the $p^2$-fold cover $F$ of ${\cal B}_K=S^{2}(p,p,p)$,
by a composition  of two $p$-fold cyclic covers.

\begin{figure}
\begin{center}
\includegraphics{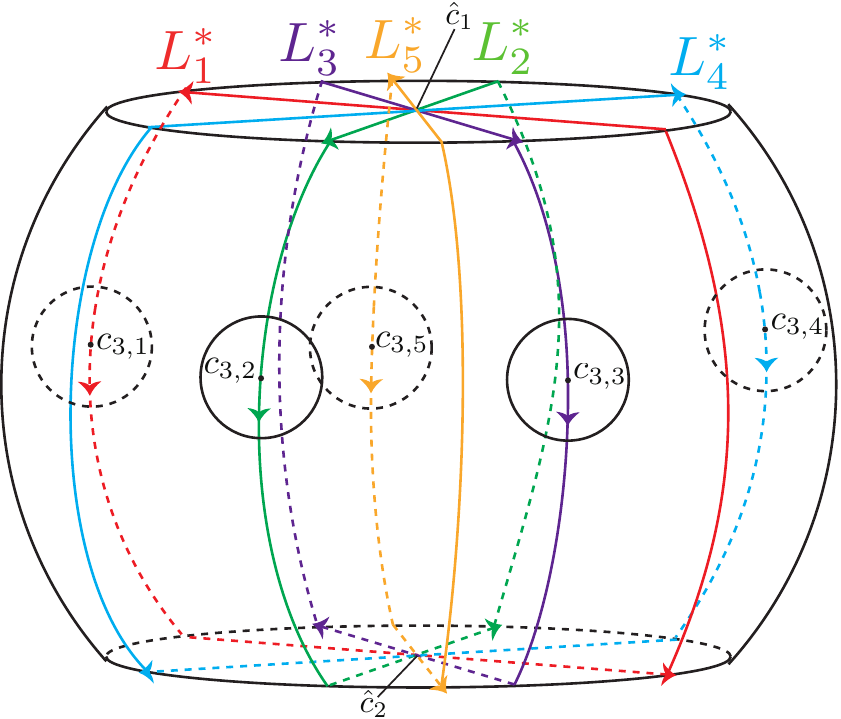}
\end{center}
\caption{\label{csF}Covering space $F'$}
\end{figure}

Let $\Gamma_1=\pi_1(S^2(p, p, p))$ be the orbifold fundametal group
of $S^2(p,p,p)$. It has a presentation
\begin{equation*}
\Gamma_1=<x_1, x_2, x_3\ |\ x_1^p=x_2^p=x_3^p=x_1x_2x_3=1>.
\end{equation*}
where $x_r$ is represented by a small circular loop in $S^2(p, p, p)$
centered at $c_r, r=1, 2, 3$. Let $\psi_1: F'\rightarrow S^2(p, p, p)$ be
the $p$-fold cyclic orbifold cover of $S^2(p, p, p)$ corresponding to the homomorphism:
\begin{equation*}
h_1: \Gamma_1\rightarrow \mathbb{Z}/p
\end{equation*}
where $h_1(x_1)=\overline{1}, h_1(x_2)=-\overline{1}$, and
$h_1(x_3)=\overline{0}$. $F'$ is as shown in Figure \ref{csF} (while $p=5$).

The order of $h_1(x_1)$ and $h_1(x_2)$ are both $p$, so
$\psi_1^{-1}(c_1)=\hat{c}_1$ and $\psi_1^{-1}(c_2)=\hat{c}_2$ are
two points in $F'$ (not cone points). The order of $h_1(x_3)$ is
$0$, so $\psi_1^{-1}(c_3)=\{c_{3, 1}, c_{3, 2},\cdots, c_{3, p}\}$
are the only cone points in $F'$, each has order $p$. We have
$F'=S^2(\underbrace{p, p, \cdots, p)}_{p}$.
 $\psi_1^{-1}(K^*)$ is a set of $p$ geodesics $L_1^*, L_2^*,\cdots, L_p^*$,
  such that $L_i^*$  goes through the cone point $c_{3,i}$ (cf. Figure \ref{csF}).

Let $\tau_1^{*}$ be the deck transformation of $\psi_1$ corresponding to $\bar{1}\in \mathbb{Z}_p$. Fix$(\tau_1^*)=\{\hat{c}_1, \hat{c}_2\}$. We may assume that $\tau_1^*(L_i^*)=L_{i+1}^*$. Orient $L_1^*$ and give $L_i^*$ the induced orientation, $1<i\leqslant p$. $F'$ admits an orientation such that $\tau_1^*$ is a counterclockwise rotation by ${2\pi}/{p}$ near $\hat{c}_1$ and a clockwise rotation by ${2\pi}/{p}$ near $\hat{c}_2$ on $F'$, since $h_1(x_1)=\bar{1}$ and $h_1(x_2)=-\overline{1}$. Note that
$L_{i}^*$ intersects $L_j^*$ at $\hat{c}_1$ in an angle of ${2\pi(i-j)}/{p}$, also intersects $L_j^*$ at $\hat{c}_2$ in an angle of ${2\pi(j-i)}/{p}$. (We always suppose that the counterclockwise direction is the positive direction of angles.)

Let $\Gamma_2$ be the orbifold  fundamental group of $F'$. It has
the presentation
\begin{equation*}
\Gamma_2=\pi_1(S^2(\underbrace{p, p,\cdots, p}_{p}))=<y_1, y_2, \cdots, y_p\ |\ y_1^p=y_2^p= \cdots =y_p^p=y_1y_2\cdots y_p=1>,
\end{equation*}
where $y_j$ is represented by a small circular loop in $S^2(\underbrace{p, p, \cdots, p}_{p})$ centered at $c_{3, j}$. Let $\psi_2: F\rightarrow F'$ be the $p$-fold cyclic orbifold cover of $F'$ corresponding to the homomorphism:
\begin{equation*}
h_2:\Gamma_2 \rightarrow \mathbb{Z}/p
\end{equation*}
where $h_2(y_j)=\overline{1}$, $1\leqslant j\leqslant p$.

The order of $h_2(y_j)$ is $p$, so
$\psi_2^{-1}(c_{3,j})=\hat{c}_{3,j}$ is a point, $1\leqslant
j\leqslant p$.  Hence  $F$ is a smooth closed orientable surface
without cone points. For each $r=1,2$, $\psi_2^{-1}(\hat{c}_r)$
is a set of $p$ points, which we denote by  $\hat{c}_{r, 1},
\hat{c}_{r, 2}, \cdots, \hat{c}_{r, p}$.  For each $i=1,...,p$,
$\psi_2^{-1}(L_i^*)$ is a set of $p$ simple closed geodesics, which
we denote by $L_{i, 1}^*, L_{i, 2}^*, \cdots, L_{i, p}^*$.

Denote the deck transformation of $\psi_2$ corresponding to $\bar{1}\in \mathbb{Z}/p$
 by $\tau_2^*$. Fix$(\tau_2^*)=\{\hat{c}_{3, 1}, \hat{c}_{3, 2},\cdots, \hat{c}_{3, p}\}$.
 We may assume that $\tau_2^*(L_{i, j}^*)=L_{i, j+1}^*$, $\tau_2^*(\hat{c}_{1, j})=\hat{c}_{1, j+1}$, and $\tau_2^*(\hat{c}_{2, j})=\hat{c}_{2, j+1}$.
Then $\psi_2^{-1}(L_i^*)=\{L_{i,j}^*:1\leqslant j\leqslant p\}$.
Now we fix an orientation for each of $L_{1, 1}^*, L_{2, 2}^*, \cdots, L_{p, p}^*$  such that $L_{j,j}^*$ goes through $\hat{c}_{3,j}$, $\hat{c}_{2,j}$ and $\hat{c}_{1,j}$ in order, and give $L_{i, j}^*=(\tau_2^*)^{j-i}(L_{i, i}^*)$ the induced orientation.
$F$ admits an orientation such that $\tau_2^*$ is a counterclockwise rotation by $2\pi/p$ near the fixed points $\hat{c}_{3, i}$, so $L_{i,1}^*, L_{i,2}^*, \cdots, L_{i,p}^*$ intersect at $\hat{c}_{3,i}$, and $L_{i,j}^{*}$ intersects $L_{i,i}^{*}$ at $\hat{c}_{3,i}$ in an angle of $2\pi(j-i)/p$, $1\leqslant i,j\leqslant p$.

\begin{figure}
\begin{center}
\includegraphics[width=6in]{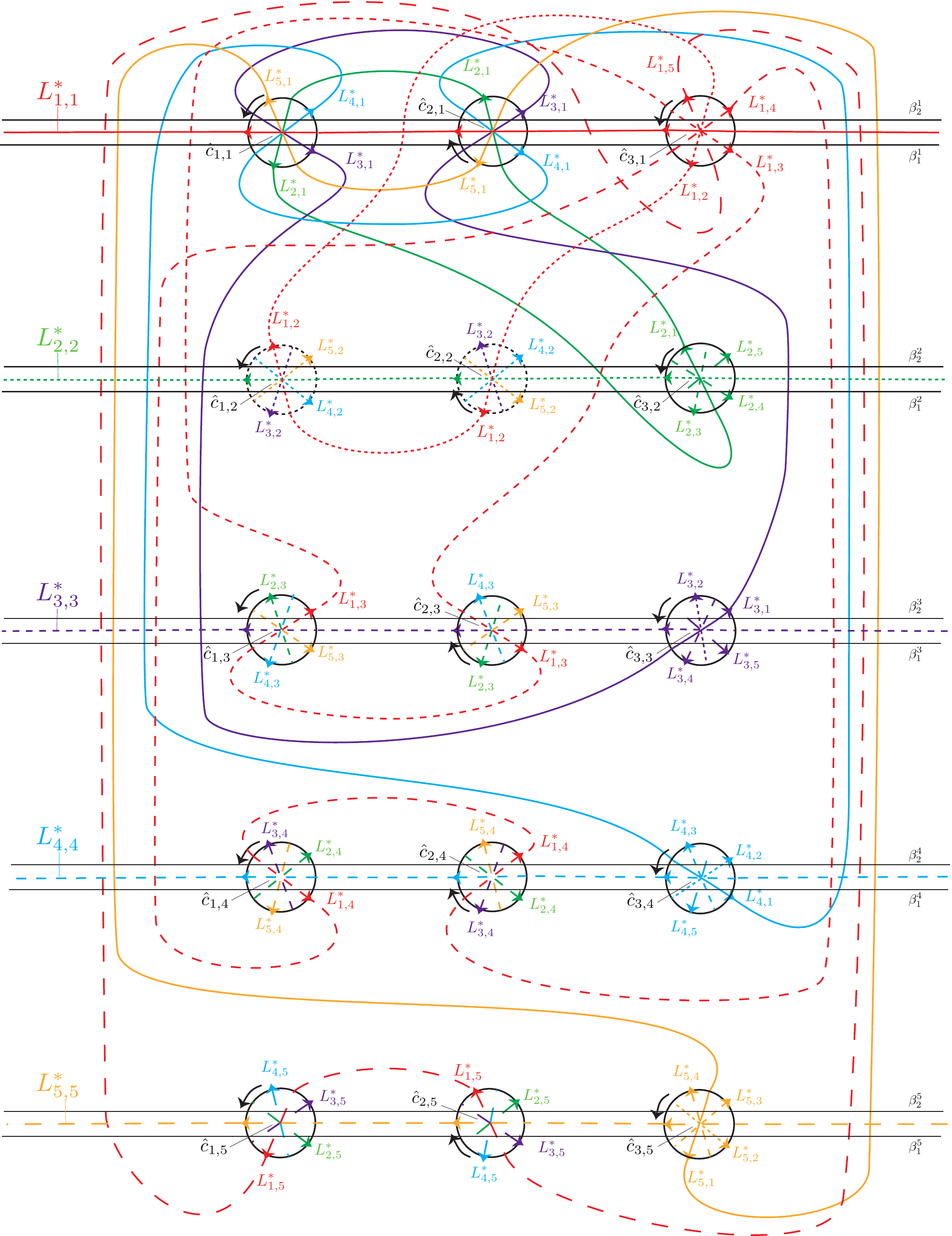}
\end{center}
\caption{\label{Lij*}$L_{i,j}^*,1\leqslant i,j \leqslant p$ with $p=5$}
\end{figure}

Since $\hat{c}_r=\psi_2(\hat{c}_{r,j})$ is not a cone point on $F'$, a small regular neighborhood of $\hat{c}_{r, j}$ on $F$ is a copy of a small regular neighborhood of $\hat{c}_r$ on $F'$, $r=1, 2$.
Recall that $L_i^*$ intersects $L_j^*$ at $\hat{c}_1$ in an angle of $2\pi(i-j)/p$, and at $\hat{c}_2$ in an angle of $2\pi(j-i)/p$.
Then $L_{i,j}^*$ intersects $L_{j,j}^*$ at $\hat{c}_{1,j}$ in an angle of $2\pi(i-j)/p$, and at $\hat{c}_{2,j}$ in an angle of $2\pi(j-i)/p$.
See the schematic picture, Figure \ref{Lij*}, for example $p=5$. Here we didn't depict the genus of the surface, so some curves meet in the picture actually do not meet. For convenience, we only draw a part of $L_{j,j}^*$ in Figure \ref{Lij*}, $1\leqslant j\leqslant p$.

Summarizing  the above discussion, we have the following remark.

\begin{remark}\label{Lij}
$\{L_{i,j}^*: 1\leqslant i,j \leqslant p\}$ only intersects at the
points $\{\hat{c}_{r,j}, r=1,2,3, 1\leqslant j\leqslant p\}$ on $F$,
and $L_{i,j}^*$ goes through $\hat{c}_{3,i}$, $\hat{c}_{2,j}$ and
$\hat{c}_{1,j}$ in order following the given orientation of
$L_{i,j}^*$, $1\leqslant i,j \leqslant p$. In particular $L_{j,j}^*$
goes through $\hat{c}_{3,j}, \hat{c}_{2,j}$, and $\hat{c}_{1,j}$, so
$\{L_{j,j}^*,1\leqslant j\leqslant p\}$ are mutually disjoint, and
\begin{align*}
L_{i, j}^* \ \text{intersects} \ L_{j, j}^* \ \text{at}
\begin{cases}
\hat{c}_{1,j}  \ \text{in an angle of}\  {2\pi(i-j)}/{p}, \\
\hat{c}_{2, j} \ \text{in an angle of}\ -{2\pi(i-j)}/{p},
\end{cases}\\
L_{j, k}^*\  intersects \ L_{j, j}^* \ at \ \hat{c}_{3, j}\  \text{in an angle of}\ {-2\pi(j-k)}/{p},
\end{align*}
Note that $i,j,k,i-j,k-j$ are considered as integers mod $p$.
\end{remark}

Let $\psi=\psi_1\circ \psi_2:
F\xrightarrow{\psi_2}S^2(\underbrace{p, p, \cdots,
p}_{p})\xrightarrow{\psi_1}S^2(p, p, p)$, which is a $p^2$-fold
orbifold covering. We have $\psi^{-1}(c_r)=\{\hat{c}_{r, 1},
\hat{c}_{r, 2}, \cdots, \hat{c}_{r, p}\}$,  $r=1, 2, 3$, and
$L^*=\psi^{-1}(K^*)=\{L_{i,j}^*:1\leqslant i,j \leqslant p\}$.

Let $\Psi_1: Y' \rightarrow W_K$ be the $p$-fold cover of $W_K$
corresponding to $\psi_1$, $\Psi_2: Y\rightarrow Y'$ the $p$-fold
cover of $Y'$ corresponding to $\psi_2$, and $\Psi=\Psi_1\circ
\Psi_2$. Then $\Psi: Y\rightarrow W_K$ is a $p^2$-fold cover of
$W_K$ corresponding to $\psi$. $Y$ has locally-trivial circle bundle
Seifert structure with base surface $F$.  Set
$L=\Psi^{-1}(\widetilde{K})$, which is a geodesic link in $Y$. We have
the following commutative diagram:

\begin{center}
\includegraphics{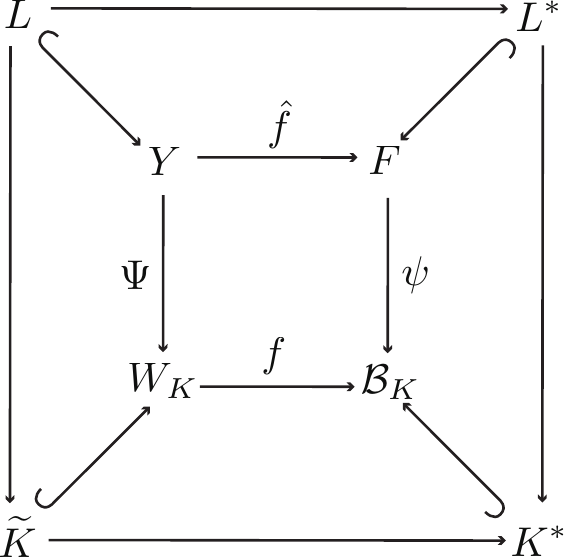}
\end{center}

Note that $L$ has exactly $p^2$ components. Let $L=\{L_{i, j}:
1\leqslant i,j \leqslant p\}$, where $L_{i, j}=\hat{f}^{-1}(L_{i,
j}^*), 1\leqslant i, j \leqslant p$. $\hat{f}|: L_{i, j} \rightarrow
L_{i, j}^*$ is a 2-fold cover.

\noindent\textbf{Step 2}

By Remark \ref{Lij}, $\{L_{j,j}^*,1\leqslant j\leqslant p\}$ are
mutually disjoint. Let $F_2^j=L_{j,j}^*\times [-\epsilon,
\epsilon]$, for some small positive number $\epsilon$,
  be a neighborhood of $L_{j,j}^*$ in $F$, such that $F_2^j$'s are mutually disjoint.
  Let  $\beta_1^j=L_{j,j}^*\times\{-\epsilon\}$,
$\beta_2^j=L_{j,j}^*\times\{\epsilon\}$, and
$F_1=\overline{F-(\underset {j=1}{\overset{p}{\cup}}F_2^j)}$.
We may suppose that $\beta_1^j$ is on the left side of $L_{j,j}^*$.
 By Remark \ref{Lij}, $L_{i,k}^*$ goes through $\hat{c}_{1.k}$,
$\hat{c}_{2,k}$ and $\hat{c}_{3,i}$, so $L_{i,k}^*$ is separated
into $2n=6$ arcs by $F_1, F_2^k$, and $F_2^i$. We denote the six
arcs of $L_{i,k}^*$ by $(L_{i,k}^1)^*$, $(L_{i,k}^2)^*$, $\cdots$,
$(L_{i,k}^6)^*$, so that $(L_{i,k}^2)^*\cup(L_{i,k}^4)^*\subset
F_2^k$, $(L_{i,k}^6)^*\subset F_2^i$, $(L_{i,k}^1)^*\cup
(L_{i,k}^3)^*\cup (L_{i,k}^5)^*\subset F_1$, and $\hat{c}_{1,k}\in
(L_{i,k}^2)^*$, $\hat{c}_{2,k}\in (L_{i,k}^4)^*$, and
$\hat{c}_{3,i}\in (L_{i,k}^6)^*$. We also require that as we travel
along $L_{i,k}^*$ in its orientation, we will pass $(L_{i,k}^6)^*$,
$(L_{i,k}^5)^*$, $\cdots$, $(L_{i,k}^1)^*$ consecutively. From the
construction, we have Remark \ref{arcs} and Table \ref{L*}.

\begin{remark}\label{arcs}
$\begin{cases}
  \text{$(L_{i,k}^1)^*$ is the part of $L_{i,k}^*$ between $\hat{c}_{3,i}$ and $\hat{c}_{1,k }$ inside $F_1$}; \\
   \text{$(L_{i,k}^3)^*$ is the part of $L_{i,k}^*$ between $\hat{c}_{1,k}$ and $\hat{c}_{2,k}$ inside $F_1$};\\
  \text{$(L_{i,k}^5)^*$ is the part of $L_{i,k}^*$ between $\hat{c}_{2,k}$ and $\hat{c}_{3,i}$ inside $F_1$}.
\end{cases}$
\end{remark}

\begin{table}
\begin{center}
\newcommand{\rb}[1]{\raisebox{-1.5ex}[0pt]{#1}}
\begin{tabular}{|c|c|c||c|c|c||c|c|c|}
\hline
$\rb{$(L_{i,k}^1)^*$}$ & $\rb{tail}$ &$\rb{head}$ & $\rb{$(L_{i,k}^3)^*$}$ & $\rb{tail}$ & $\rb{head}$& $\rb{$(L_{i,k}^5)^*$}$& $\rb{tail}$ & $\rb{head}$\\[3ex]
\hline
 $\rb{$i-k\leqslant\frac{p-1}{2}$}$ & $\rb{$\beta_1^k$}$ & $\rb{$\beta_1^i$}$ & $\rb{$i-k\leqslant\frac{p-1}{2}$}$ & $\rb{$\beta_2^k$}$ & $\rb{$\beta_2^k$}$ & $\rb{$i-k\leqslant\frac{p-1}{2}$}$ & $\rb{$\beta_2^i$}$ & $\rb{$\beta_1^k$}$ \\[3ex]
\hline
 $\rb{$i-k>\frac{p-1}{2}$}$ & $\rb{$\beta_2^k$}$ & $\rb{$\beta_2^i$}$ &  $\rb{$i-k>\frac{p-1}{2}$}$&$\rb{$\beta_1^k$}$ & $\rb{$\beta_1^k$}$ &  $\rb{$i-k>\frac{p-1}{2}$}$&$\rb{$\beta_1^i$}$ & $\rb{$\beta_2^k$}$ \\[3ex]
\hline
\end{tabular}
\end{center}
Note: $i-k$ is considered as a \textit{nonnegative} integer mod $p$\\
\caption{\label{L*}The induced orientation on $(L_{i, k}^{2l-1})^*, 1\leqslant l\leqslant 3$.}
\end{table}

We need Table \ref{L*} in the following lemma and in Step 3.

\begin{lemma}\label{connected}
$F_1$ is connected.
\end{lemma}

\noindent\textbf{Proof}: It suffices to prove that the boundary components of
$F_1$, which is the set $\{\beta_1^i, \beta_2^j, 1\leqslant i, j
\leqslant p\}$, can be mutually connected to each other by arcs in
$F_1$.

 From Table \ref{L*}, the arc $(L_{j+1,j}^5)^*$ connects
$\beta_2^{j+1}$ and $\beta_1^{j}$, and the arc $(L_{j+1,j}^1)^*$
connects $\beta_1^{j}$ and $\beta_1^{j+1}$, since
$j+1-j=1\leqslant\displaystyle{\frac{p-1}{2}}$, where $p\geq3$.
So $\beta_2^{j+1}$ and $\beta_1^{j+1}$ can be connected in $F_1$,
$1\leqslant {j+1}\leqslant p$. Also $\beta_2^{j+1}$ and $\beta_1^j$
are connected by the arc $(L_{j+1, j}^5)^*$. Hence  $\{\beta_1^i,
\beta_2^j, 1\leqslant i, j \leqslant p\}$ can be mutually connected
to each other in $F_1$. $\square$

\begin{figure}
\begin{center}
\includegraphics[width=2.8in]{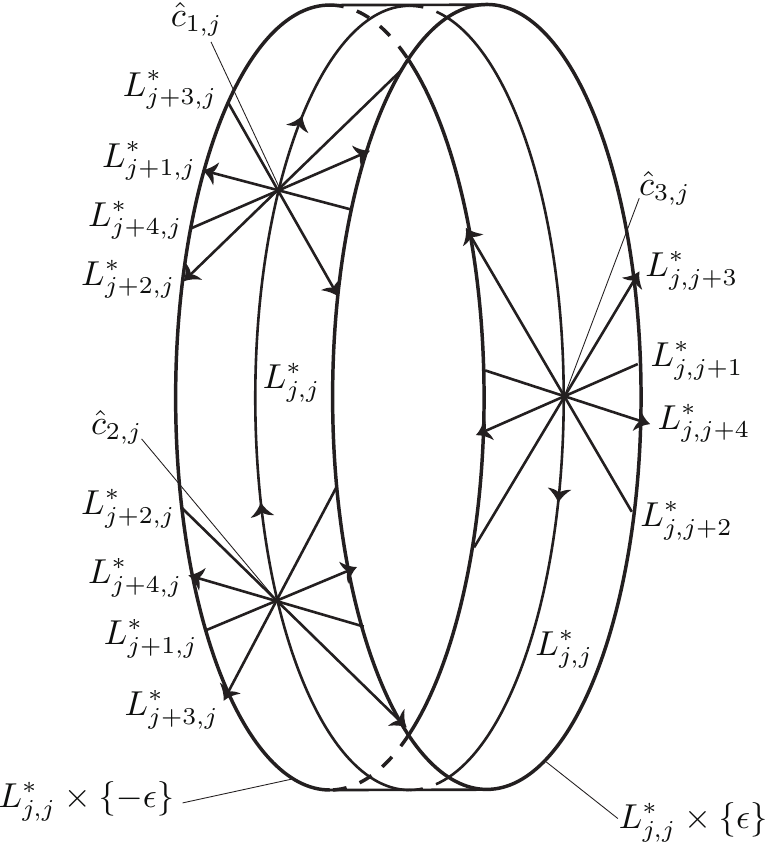}
\end{center}
\begin{center}
\caption{\label{nbhd} A neighborhood of $L_{j, j}^*$ in $F$.}
\end{center}
\end{figure}

Let $T^j=\hat{f}^{-1}(L_{j, j}^*)\subset Y$, which is a vertical
torus over $L_{j,j}^*$. Then $L_{i, j}$ and $L_{j, i}$ are
transverse to $T^j$ for all $1\leqslant i\leqslant p$ and $i\neq j$.
The situation near $L_{j, j}^*$ in $F$ is described in Figure \ref{nbhd},
while $p=5$.

Since $L_{j,j}^*$ is a geodesic, the torus  $T^j$ is a totally geodesic torus which inherits a Euclidean structure from the $\widetilde{SL_2}$ structure on $Y$. For any two simple closed geodesics $\{a_j, b_j\}\subset T^j$ with  $H_1(T^j)=<a_j,b_j>$, $T^j$ can be identified to $S^1\times S^1$ where each
$S^1\times \{*\}$ is a geodesic isotopic to $a_j$ and each $\{*\}\times S^1$ is a geodesic isotopic to $b_j, 1\leqslant j\leqslant p$.

Similar to Proposition 6.1 in \cite{abz}, we have the following
proposition.

\begin{prop}\label{semibundle}
The exterior of $\{ L_{1, 1} \cup L_{2, 2} \cup \cdots \cup L_{p, p} \}$ in
$Y$ is a surface semi-bundle.
\end{prop}

\noindent\textbf{Proof}: To prove this, we apply \cite{wy} and  compare \cite{abz}.

Let $Y_2^j$ be the submanifold of $Y$ lying over
$F_2^j=L_{j, j}^*\times [-\epsilon, \epsilon]$, and
$Y_1$ be the submanifold of $Y$ lying over $F_1$. Note
that $Y_1$ is connected since $F_1$ is connected by Lemma \ref{connected}.
Define
\begin{equation*}
T_1^j=\hat{f}^{-1}(\beta_1^j), \ T_2^j=\hat{f}^{-1}(\beta_2^j).
\end{equation*}

Let $Y_0$ be the 3-manifold obtained by cutting $Y$ open along
$T_2^j$, and $F_0$ the surface obtained by cutting $F$ open along
all $\beta_2^j$,  $1\leqslant j\leqslant p$. The restriction of the
Seifert fiberation of $Y$ to each of  $Y_0, Y_1, Y_2^j, j=1,...,p$,
is a trivial circle bundle. We give the circle fibers of $Y_0$ a
consistent orientation. Choose a horizontal section $B_0$ of the
bundle $Y_0\rightarrow F_0$ such that $T^j\cap B_0$ is a geodesic,
$1\leqslant j\leqslant p$. Let $B_1$, $B_2^j$ be the restriction of
$B_0$ in $Y_1$ and $Y_2^j$ respectively. Fix an orientation of $B_0$
and let $B_1, B_2^j$ have the induced orientation. And let $\partial B_1$
and $\partial B_2^j$ have the induced orientation.

We denote the torus $T_r^j $ by  $T_{1, r}^j$ and $T_{2, r}^j$ when
we think of it as lying in $Y_1$ and $Y_2^j$ respectively, $r=1, 2$. Let $\phi
_{k, r}^j$ be a fixed circle fiber  in the torus $T_{k,r}^j$ for
each of $r=1, 2; k=1, 2; 1\leqslant j\leqslant p$. Let
$\alpha_{k,r}^j=T_{k,r}^j\cap B_0$. Then  $\alpha _{k, 1}^j=-\alpha
_{k, 2}^j$, $\phi _{k, 1}^j=\phi _{k, 2}^j$, and
$\{\alpha _{k, r}^j, \phi _{k, r}^j\}$ form a basis of $H_1(T_{k,
r}^j)$, $k=1, 2$; $r=1, 2$; $1\leqslant j\leqslant p$.

By choosing the proper horizontal section $B_0$, we can assume that the Seifert manifold $Y$ is obtained form $Y_1$ and $Y_2^j$'s by
gluing $T_{1, 1}^j$ to $T_{2, 1}^j$ and $T_{1, 2}^j$ to $T_{2, 2}^j$
using maps $g_r^j: T_{2, r}^j\rightarrow T_{1, r}^j\ (r= 1, 2)$ determined by the conditions
\begin{equation*}
(g_1^j)_\ast(\alpha _{2, 1}^j)=-\alpha _{1, 1}^j\ \ \ \ \ \ (g_1^j)_\ast(\phi _{2, 1}^j)=\phi _{1, 1}^j
\end{equation*}
\begin{equation*}
(g_2^j)_\ast(\alpha _{2, 2}^j)=-\alpha _{1, 2}^j+e^j\phi _{1, 2}^j\ \ \ \ \ \ (g_2^j)_\ast(\phi _{2, 2}^j)=\phi _{1, 2}^j
\end{equation*}
where $e= \sum_{j=1}^pe^j\in \mathbb{Z}$ is the Euler number of the oriented circle bundle $Y\rightarrow F$. $Y\rightarrow W_K$ is a $p^2$-fold cover so
\begin{equation*}
e=p^2e(W_K)=-p^2(\frac{q_1}{p}+\frac{q_2}{p}+\frac{q_3}{p})=-p(q_1+q_2+q_3).
\end{equation*}
Note that $q_1+q_2+q_3$ is an odd number since $K$ has only one
component. For convenience, we may assume $e^1=e^2=\cdots
=e^p=\displaystyle{\frac{e}{p}}=-(q_1+q_2+q_3)=\tilde{e}$, so
$\tilde{e}$ is odd.

Each $Y_2^j, j=1,...,p$ also has a circle fibration with $L_{j,j}$
as a fiber. As in \cite{abz}, we call the circle fibers of this circle fibration {\it new fibers} and called the fibers of $Y$ {\it
original fibers}. Give new fibers of $Y_2^j$ a fixed consistent
orientation.
 We  shall denote by $\overline{\phi}^j$ a general new fiber in $M_2^j$, and
$\phi$ a general  original  fiber of $Y$.

Let $N^j$ be the regular neighborhood of $L_{j, j}$ in $Y_2^j$,
which is disjoint from other components of $L$ and consists of new
fibers of $Y_2^j$. Let $M_2^j=Y_2^j\setminus\overset{\circ}{N^j}$.
Then $M_2^j$ is a Seifert fibered space whose circle fibers are new
fibers of $Y_2^j$. Set  $T_3^j=\partial {N^j}$. The exterior of
$\{L_{1, 1}, L_{2, 2},\cdots , L_{p, p}\}$ in $Y$, $M$, is a graph
manifold and has the following JSJ decomposition
\begin{equation*} M=Y_1\cup M_2^1\cup M_2^2\cup \cdots
\cup M_2^p
\end{equation*}

Let $\overline{B}_2^j$ be the image of one section of the circle fibration
of $M_2^j$, such that $\overline{B}_2^j$ intersects $T^j$ in a
geodesic. Fix an orientation for $\overline{B}_2^j$ and let $\partial
\overline{B}_2^j$ have the induced orientation.
 There is another
basis of $H_1(T_{2, r}^j)$, $\{\overline{\alpha}_{2, r}^j,
\overline{\phi}_{2, r}^j\}$, where $\overline{\phi}_{2, r}^j$ is a
fixed new fiber on $T_{2,r}^j$  and $\overline{\alpha}_{2, r}^j$ is
the component of $\partial \overline{B}_2^j$ on $T_{2, r}^j$, $r=1,
2$, with their chosen orientations.
 Let $\overline{\alpha}_{2, 3}^j=\overline{B}_2^j\cap T_3^j$.

 The relation between the old basis
 $\{\alpha_{2, 1}^j, \phi_{2, 1}^j\}$ of $H_1(T_{2, 1}^j)$ and
 the new one $\{\overline{\alpha}_{2, 1}^j, \overline{\phi}_{2, 1}^j\}$
is given by
\begin{equation*}
\overline{\alpha}_{2, 1}^j=a^j\alpha_{2, 1}^j+b^j\phi_{2, 1}^j\ \ \ \ \ \ \overline{\phi}_{2, 1}^j=c^j\alpha_{2, 1}^j+d^j\phi_{2, 1}^j
\end{equation*}
where $a^j, b^j, c^j, d^j$ are integers satisfying $a^jd^j-b^jc^j=\pm 1$. We may assume that $a^jd^j-b^jc^j=1$ by reversing the orientation of the fibres $\phi$ if necessary. For convenience, we can suppose that $a^1=a^2= \cdots =a^p=a,\ b^1=b^2=\cdots =b^p=b,\ c^1=c^2=\cdots =c^p=c$, and $d^1=d^2=\cdots =d^p=d$.

In $M_2^j$, we have $\overline{\alpha}_{2, 1}^j=-\overline{\alpha}_{2, 2}^j$ and $\overline{\phi}_{2, 1}^j=\overline{\phi}_{2, 2}^j$. Thus
\begin{equation*}
\overline{\alpha}_{2, 2}^j=a\alpha_{2, 2}^j-b\phi_{2, 2}^j\ \ \ \ \ \ \overline{\phi}_{2, 2}^j=-c\alpha_{2, 2}^j+d\phi_{2, 2}^j.
\end{equation*}

Hence, with respect to the basis $\{\alpha_{1, r}^j, \phi_{1, r}^j\}$ of $H_1(T_{1, r}^j)$ and $\{\overline{\alpha}_{2, r}^j, \overline{\phi}_{2, r}^j\}$ of $H_1(T_{2, r}^j)\ (r=1, 2)$, the gluing maps $g_1^j: T_{2, 1}^j\rightarrow T_{1, 1}^j$ and $g_2^j: T_{2, 2}^j\rightarrow T_{1, 2}^j$ can be expressed as
\begin{equation*}
(g_1^j)_*(\overline{\alpha}_{2, 1}^j)=-a\alpha_{1, 1}^j+b\phi_{1, 1}^j\ \ \ \ \ \ (g_1^j)_*(\overline{\phi}_{2, 1}^j)=-c\alpha_{1, 1}^j+d\phi_{1, 1}^j
\end{equation*}
\begin{equation*}
(g_2^j)_*(\overline{\alpha}_{2, 2}^j)=-a\alpha_{1, 2}^j+(a\widetilde{e}-b)\phi_{1, 2}^j\ \ \ \ \ \ (g_2^j)_*(\overline{\phi}_{2, 2}^j)=c\alpha_{1, 2}^j+(d-c\widetilde{e})\phi_{1, 2}^j
\end{equation*}
The associated matrices are
\begin{equation*}
G_1^j=(g_1^j)_*=\begin{pmatrix}-a&b\\-c&d\end{pmatrix}\ \ \ \ \ \ G_2^j=(g_2^j)_*=\begin{pmatrix}-a&{a\widetilde{e}-b}\\c&{d-c\widetilde{e}}\end{pmatrix}
\end{equation*}
\begin{equation*}
(G_1^j)^{-1}=\begin{pmatrix}-d&b\\-c&a\end{pmatrix}\ \ \ \ \ \ (G_2^j)^{-1}=\begin{pmatrix}{c\widetilde{}e-d}&{a\widetilde{e}-b}\\c&a\end{pmatrix}.
\end{equation*}

\begin{figure}
\begin{center}
\includegraphics{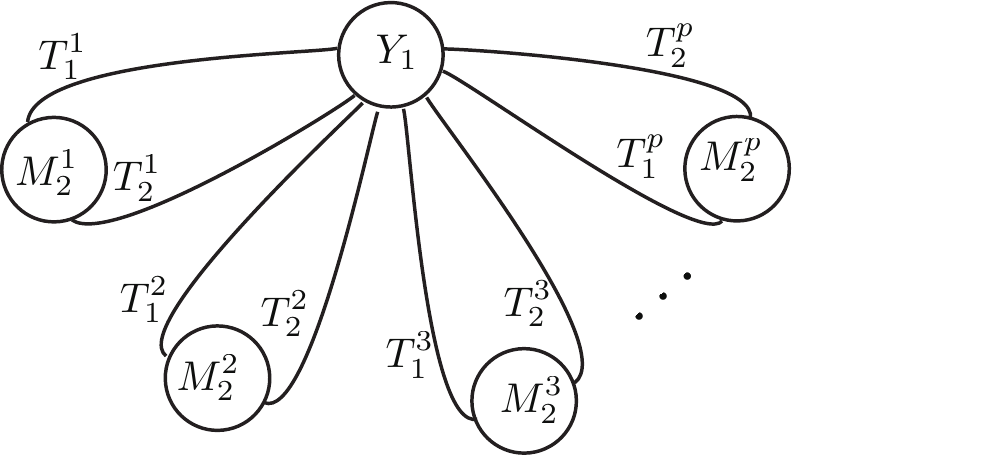}
\end{center}
\caption{\label{JSJ}The graph decomposition of $M$.}
\end{figure}

 The graph of the JSJ-decomposition of $M$ consists of
 $p+1$ vertices corresponding to
 $Y_1, M_2^1, M_2^2, \cdots, M_2^p$, and $2p$ edges
  corresponding to $T_1^j$'s and $T_2^j$'s. See Figure \ref{JSJ}.%

By \cite{wy}, we can get a horizontal surface of $M$ from every non-zero
solution of the following equation.
\begin{equation*}
(Y-Z)\left(\begin{array}{c}\lambda \\ \bar{\lambda}_1 \\  \vdots \\ \bar{\lambda}_p\end{array}\right)=\left(\begin{array}{c}0\\{\ast}\\ \vdots \\{\ast} \end{array}\right),\ \lambda, \bar{\lambda}_1, \cdots, \bar{\lambda}_p, \ast \in \mathbb{Z}
\end{equation*}
where  $Y$ and $Z$ are $(p+1)\times(p+1)$ matrices defined on
 page 450 of \cite{wy}.
  The entries of $Y$ and $Z$ are decided by the gluing matrix.
From \cite{wy}, we have $y_{1, i}=y_{i,
1}=\displaystyle{\frac{2}{c}}, 2\leqslant i\leqslant {p+1}$, and
other entries of $Y$ are all zeroes.
$Z=diag(z_1,z_2,\cdots,z_{p+1})$ is a diagonal matrix with
$z_1=\sum_{i=1}^p(\displaystyle{\frac{d}{c}+\frac{c\widetilde{e}-d}{c}})=p\widetilde{e}=e,
z_i=\displaystyle{\frac{a}{c}-\frac{a}{c}}=0, 2\leqslant i\leqslant
{p+1}$. The matrix equation (1.6)
of \cite{wy} becomes:

\begin{equation*}
\begin{pmatrix}{-e}&{\displaystyle{\frac{2}{c}}}&{\displaystyle{\frac{2}{c}}}&{\cdots}&{\displaystyle{\frac{2}{c}}}\\
{\displaystyle{\frac{2}{c}}}&0&0&{\cdots}&0\\[2ex]
{\displaystyle{\frac{2}{c}}}&0&0&{\cdots}&0\\{\vdots}&{\vdots}&{\vdots}&{\vdots}&{\vdots}\\{\displaystyle{\frac{2}{c}}}&0&0&{\cdots}&0\end{pmatrix}_{(p+1)\times (p+1)}\begin{pmatrix}\lambda\\\overline{\lambda}_1\\\overline{\lambda}_2\\{\vdots}\\\overline{\lambda}_p\end{pmatrix}=\begin{pmatrix}0\\{\ast}\\{\ast}\\{\vdots}\\{\ast}\end{pmatrix}\Rightarrow \begin{cases}{-e\lambda +\displaystyle{\frac{2}{c}}\overline{\lambda}_1+\displaystyle{\frac{2}{c}}\overline{\lambda}_2+{\cdots}+\displaystyle{\frac{2}{c}}\overline{\lambda}_p=0}\\{\displaystyle{\frac{2}{c}}\lambda=\ast}\end{cases}
\end{equation*}

We may assume $\overline{\lambda}_j=\overline{\lambda}$, $1\leqslant j\leqslant p$. Then

\begin{equation}\label{es}
\begin{cases}{-e\lambda+\displaystyle{\frac{2p}{c}}\overline{\lambda}=0}\\{\displaystyle{\frac{2}{c}}\lambda=\ast}
\end{cases}
\end{equation}

We can construct horizontal surfaces $H_1, H_2^j$ from the solution of (\ref{es}) such that the projection of $H_1$ to $F_1$ has degree $|\lambda|$ and that of $H_2^j$ to the base of $M_2^j$ has degree $|\overline{\lambda}|$ (\cite{wy}). Suppose $\lambda\neq 0$. Then

\begin{equation*}
\frac{\lambda}{\bar{\lambda}}=\frac{2p}{ce}=\frac{2}{c\tilde{e}},\ \ \ \ \frac{\bar{\lambda}}{\lambda}=\frac{ce}{2p}=\frac{c\tilde{e}}{2}.
\end{equation*}

Suppose that $\partial H_1$ on $T_i^j$ with respect to the basis $\{\alpha_{1, i}^j, \phi_{1, i}^j\}$ is $u_i^j\alpha_{1, i}^j+t_i^j\phi_{1, i}^j\ (i=1, 2)$, and $\partial H_2^j$ on $T_i^j$ with respect to the basis $\{\bar{\alpha}_{2, i}^j, \bar{\phi}_{2, i}^j\}$ is $\bar{u}_i^j\bar{\alpha}_{2, i}^j+\bar{t}_i^j\bar{\phi}_{2, i}^j\ (i=1, 2, 3)$. Then there are $\epsilon_1^j,\epsilon_2^j\in\{\pm1\}$, so that equation (1.2) of \cite{wy} takes on the form:

\begin{equation*}
\frac{t_1^j}{u_1^j}=\frac{\epsilon_1^j\bar{\lambda}}{-\lambda c}-\frac{d}{c}=\frac{-\epsilon_1^jc\tilde{e}-2d}{2c},\ \ \ \ \ \ \frac{\bar{t}_1^j}{\bar{u}_1^j}=\frac{\epsilon_1^j\lambda}{-\bar{\lambda}c}+\frac{a}{-c}=\frac{-2\epsilon_1^j-ac\tilde{e}}{\tilde{e}c^2},
\end{equation*}

\begin{equation*}
\frac{t_2^j}{u_2^j}=\frac{\epsilon_2^j\bar{\lambda}}{\lambda c}-\frac{c\tilde{e}-d}{c}=\frac{2d-(2-\epsilon_2^j)c\tilde{e}}{2c},\ \ \ \ \ \ \frac{\bar{t}_2^j}{\bar{u}_2^j}=\frac{\epsilon_2^j\lambda}{\bar{\lambda}c}+\frac{a}{c}=\frac{2\epsilon_2^j+ac\tilde{e}}{\tilde{e}c^2}.
\end{equation*}

Without lose of generality, we can assume $\epsilon_1^1=\epsilon_1^2=\cdots =\epsilon_1^p=\epsilon_1, \epsilon_2^1=\epsilon_2^2=\cdots =\epsilon_2^p=\epsilon_2$.
\bigskip

We have $\sum_{j=1}^p(\displaystyle{\frac{t_1^j}{u_1^j}+\frac{t_2^j}{u_2^j}})=0$, since $H_1$ is a horizontal surface. So
\begin{equation*}
p[-\epsilon_1c\tilde{e}-2d+2d-(2-\epsilon_2)c\tilde{e}]=0
\end{equation*}
\begin{equation*}
-\epsilon_1+\epsilon_2=2.
\end{equation*}

From the above equation, $\epsilon_2=-\epsilon_1=1$. Then
\begin{equation}\label{slope12}
\parbox{3cm}
{\begin{eqnarray*}
\frac{t_1^j}{u_1^j}=\frac{c\tilde{e}-2d}{2c}, & &  \frac{\bar{t}_1^j}{\bar{u}_1^j}=\frac{2-ac\tilde{e}}{\tilde{e}c^2}; \\
\frac{t_2^j}{u_2^j}=\frac{2d-c\tilde{e}}{2c},  & & \frac{\bar{t}_2^j}{\bar{u}_2^j}=\frac{2+ac\tilde{e}}{\tilde{e}c^2}.
\end{eqnarray*}}
\hfill
\end{equation}

Also $\displaystyle{\frac{\overline{t}_1^j}{\overline{u}_1^j}+\frac{\overline{t}_2^j}{\overline{u}_2^j}+\frac{\overline{t}_3^j}{\overline{u}_3^j}}=0$, so
\begin{equation}\label{slope3}
\frac{\overline{t}_3^j}{\overline{u}_3^j}=-(\frac{\overline{t}_1^j}{\overline{u}_1^j}+\frac{\overline{t}_2^j}{\overline{u}_2^j})=-(\frac{2-ac\tilde{e}+2+ac\tilde{e}}{\tilde{e}c^2})=-\frac{4}{\tilde{e}c^2}.
\end{equation}

In fact, $L_{j, j}$ double covers $L_{j, j}^*$, so $c=2$. Then, we can take $\lambda=1, \bar{\lambda}=\tilde{e}$. As in \cite{abz}, we take $b=0, a=d=1$. Since $\tilde{e}$ is odd, we can determine $t_i^j, u_i^j, \overline{t}_i^j$ and $\overline{u}_i^j$, as following

\begin{equation}\label{t}
\parbox{10cm}
{\begin{eqnarray*}
t_1^j=\frac{\tilde{e}-1}{2},\  u_1^j=1; \;\;\;\;\;  t_2^j=\frac{1-\tilde{e}}{2}, \ u_2^j=1\;\;\;\;\;\;\;\;\;\;\;\;\;\;\;\;\;\;\\
\bar{t}_1^j=\frac{1-\tilde{e}}{2}, \ \bar{u}_1^j=\tilde{e}; \;\;\;\;\; \bar{t}_2^j=\frac{1+\tilde{e}}{2}, \ \bar{u}_2^j=\tilde{e}; \;\;\;\;\; \ \bar{t}_3^j=-1, \ \bar{u}_3^j=\tilde{e}.
\end{eqnarray*}}
\hfill
\end{equation}

$\partial H_1\cap T_{1, i}^j$ has $\displaystyle{|\frac{\lambda}{u_i^j}|=1}$ component, $i=1, 2$, and $\partial H_2^j\cap T_{2, i}^j$ has $\displaystyle{|\frac{\bar{\lambda}}{\bar{u}_i^j}|=1}$ component, $i=1, 2, 3$.  Since $\epsilon_1^j=-\epsilon_2^j$, $H=H_1\cup(\underset{j=1}{{\overset{p}{\cup}}}H_2^j)$ is non-orientable, so $M$ is a surface semi-bundle.
$\square$

Next, we give the constructions for $H_1$ and $H_2^j$ in detail.
The construction of $H_2^j$ is the same as the construction of
 $H_2$ in \cite{abz}.

Set $\overline{\alpha}_0^j=(\overline{B}_2^j\cap T^j)$. $\overline{\phi}_0^j$ is a fixed new fibre in $T^j$. $\overline{\alpha}_0^j$ and $\overline{\phi}_0^j$ are both geodesics. Then $T^j$ is identified with $\overline{\alpha}_0^j\times \overline{\phi}_0^j$ and $Y_2^j$ with $\overline{\alpha}_0^j\times \overline{\phi}_0^j\times\left[-\epsilon, \epsilon\right]$. Separate $\overline{\alpha}_0^j$ into two parts, $\overline{\alpha}_0^j=I^j\cup J^j$, where $I^j$ meets $L_{j, j}$. Take $0<\delta<\epsilon$, then $M_2^j$ can be expressed as following
\begin{equation*}
M_2^j=(\overline{\alpha}_0^j\times \overline{\phi}_0^j\times [-\epsilon, \epsilon])\setminus(\overset{\circ}{I^j}\times \overline{\phi}_0^j\times(-\delta, \delta)).
\end{equation*}

\begin{figure}
\begin{center}
\includegraphics[width=6in]{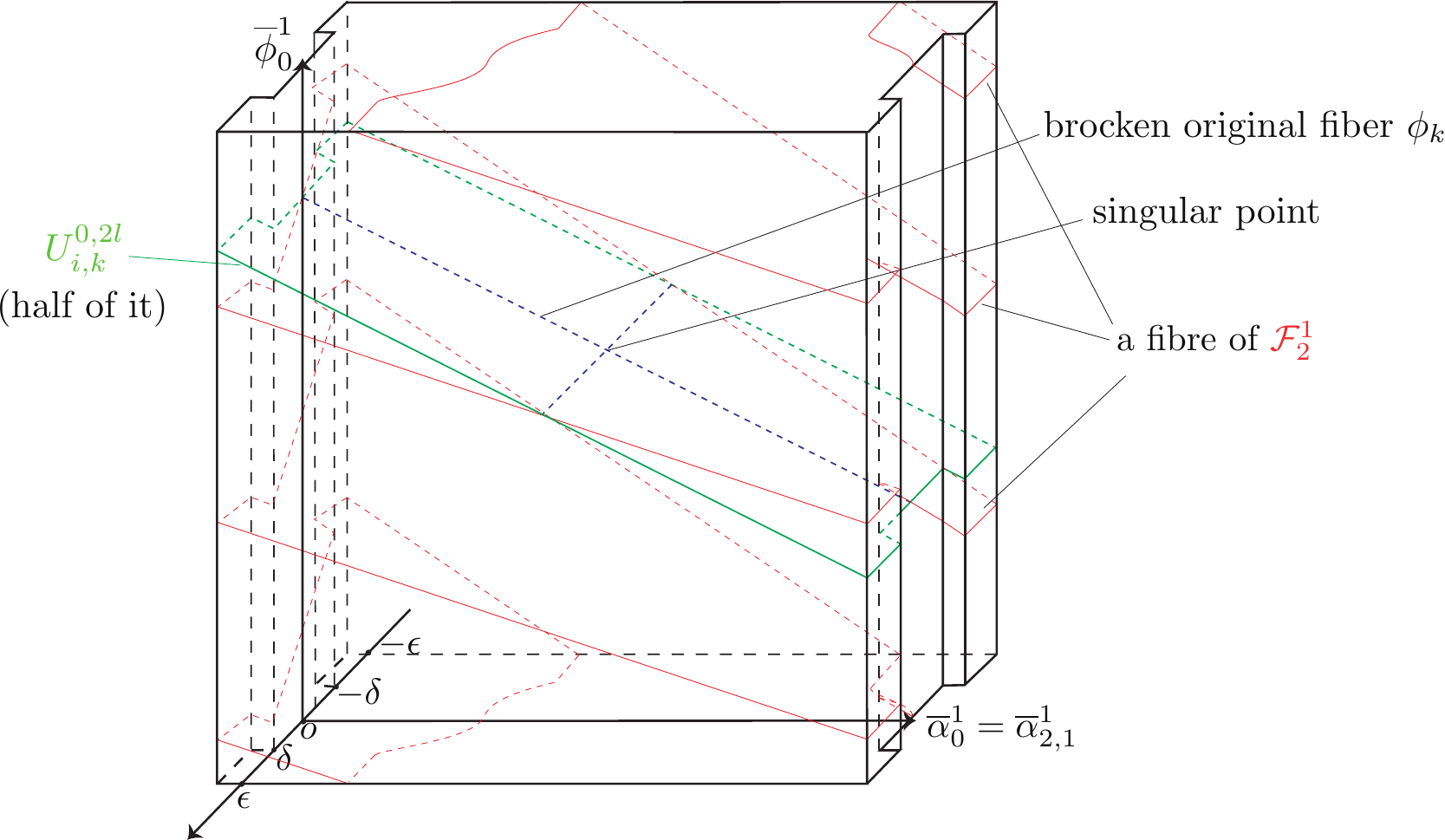}
\end{center}
\caption{\label{fibre}One fibre of $\mathcal{F}_2^1$ in $M_2^1$.}
\end{figure}

As in \cite{abz}, we have $H_2^j=\Theta_-^j\cup \Theta_0^j\cup \Theta_+^j$. $\Theta_-^j=\gamma_1\times[-\epsilon, -\delta]$, where $\gamma_1$ is a fixed connected simple closed geodesic in $T^j\times \{-\epsilon\}$ of slope $\displaystyle{\frac{1}{2\tilde{e}}-\frac{1}{2}}$ with respect to the basis $\{\overline{\alpha}_0^j\times \{-\epsilon \}, \overline{\phi}_0^j\times \{-\epsilon \}\}$. $\Theta_+^j=\gamma_2\times[\delta, \epsilon]$, where $\gamma_2$ is a fixed connected simple closed geodesic in $T^j\times \{\delta\}$ of slope $\displaystyle{-\frac{1}{2\tilde{e}}-\frac{1}{2}}$ with respect to the basis $\{\overline{\alpha}_0^j\times \{\delta \}, \overline{\phi}_0^j\times \{\delta \}\}$, see (\ref{t}). (Note that $\overline{\alpha}_0^j\times \{\delta \}=-\bar{\alpha}_{2,2}^j$.) $\Theta_0^j$ is a surface in $J^j\times \overline{\phi}_0^j\times[-\delta, \delta]$ such that for each $t\in (-\delta, \delta),\ \Theta_0^j\cap T^j\times \{t\}$ is a union of $|\tilde{e}|$ evenly spaced geodesic arcs of slope $\displaystyle{-\frac{1}{2}-\frac{t}{2\delta \tilde{e}}}$. $H_2^j$ is everywhere transverse to the new fibres. $H_2^j$ is also transverse to the original fibres except when $t=0$. Thus we can construct a surface fibration $\mathcal{F}_2^j$ in $M_2^j$ by isotoping $H_2^j$ around the new fibres. Similar as \cite{abz}, Figure \ref{fibre} illustrates one fibre of $\mathcal{F}_2^1$ in $M_2^1$, when $K=(\displaystyle{\frac{1}{5}, \frac{1}{5}, \frac{1}{5}})$, $\tilde{e}=-3$ and $\delta$ is small enough.

Next we construct $H_1$ similarly to \cite{abz}. Recall that $B_1$
is a horizontal section of $Y_1\rightarrow F_1$, and $B_1\cap (T_{1,
1}^j\cup T_{1, 2}^j)$ are geodesics. $F_1$ is connected by Lemma
\ref{connected}, so $B_1$ is connected. For each $j$, fix a properly
embedded arc $\sigma^j$ in $B_1$ connecting $T_{1, 1}^j$ and $T_{1,
2}^j$ such that all $\sigma^j$, $j=1,..,p$ are mutually disjoint
(such set of arcs obviously exists). Let $\sigma^j\times[-1, 1]$ be a
regular neighborhood of $\sigma^j$ in $B_1$ such that $(\sigma^j\cap
T_{1, 1}^j)\times[-1, 1]$ is a subarc of $\alpha_{1, 1}^j$ which
follows the orientation of $\alpha_{1, 1}^j$ as we pass from -1 to
1. Then $H_1$ is obtained by replacing every $\sigma^j\times[-1, 1]$
in $B_1$ by a reimbedding of it which wraps around the
$\phi^j$-direction $\displaystyle{(\frac{1-\tilde{e}}{2})}$ times
as we pass from -1 to 1 in $\sigma^j\times[-1, 1]\times\phi^j$. Let
 ${\cal F}_1$ be the corresponding surface fibration of $Y_1$
 with $H_1$ as a surface fiber.

 As in \cite{abz}, we may suppose that
  $\partial H_1=\partial H_2^1\cup \partial
H_2^2\cup\cdots\cup\partial H_2^p$ and that $\mathcal{F}_1\cup
\mathcal{F}_2^1\cup \cdots \cup \mathcal{F}_2^p$
  forms a semi-surface bundle $\mathcal{F}$ in $M$,
  as described in Proposition \ref{semibundle}

\noindent\textbf{Step 3}

Recall that $T^j=\hat{f}^{-1}(L_{j, j}^*)$ is vertical in the
Seifert structure of $Y$. Let $U_{i, k}=\hat{f}^{-1}(L_{i, k}^*)$,
so $L_{i, k}\subset U_{i, k}$. By Remark \ref{Lij}, we have $U_{i,
j}\cap U_{k, j}=\phi_{1, j}\cup\phi_{2, j}$, $i\neq k$, and $U_{i,
j}\cap U_{i, k}=\phi_{3, i}$, $j\neq k$, where $\phi_{r,
j}=\hat{f}^{-1}(\hat{c}_{r, j}), r=1,2$, and $\phi_{3,
i}=\hat{f}^{-1}(\hat{c}_{3, i})$. $L_{j, j}$ is a 2-fold cover of
$L_{j, j}^*$, so it intersects each of $\phi_{1, j}, \phi_{2, j}$
and $\phi_{3, j}$ twice. Then $U_{i, k}^0=U_{i, k}\cap M$ is a
6-punctured torus.

Let $U_{i,k}^{0,2l}=\hat{f}^{-1}((L_{i,k}^l)^*)\cap M$, $1\leqslant i,k\leqslant p, i\neq k, 1\leqslant l\leqslant 6$. $U_{i,k}^{0,2l-1}$ is a vertical annulus, and $U_{i,k}^{0, 2l}$ is a vertical annulus with two punctures, $l=1,2,3.$ By the way that we define $(L_{i,k}^l)^*$ in Step 2, $(U_{i, k}^{0, 1}\cup U_{i, k}^{0, 3}\cup U_{i, k}^{0, 5})\subset Y_1$, $(U_{i, k}^{0, 2}\cup U_{i, k}^{0, 4})\subset M_2^k$, $U_{i, k}^{0, 6}\subset M_2^i$ and $U_{i, k}^{0, 2l}\cap \phi_{l, k}\neq \emptyset$, $l=1, 2$, and $U_{i, k}^{0, 6}\cap \phi_{3, i}\neq \emptyset$.

Similarly, $L_{i, k}$ is separated into 12 parts, $L_{i, k}^1, L_{i, k}^2, \cdots, L_{i, k}^{12}$, such that $(L_{i, k}^l\cup L_{i, k}^{l+6})\subset U_{i, k}^{0, l}, 1\leqslant l\leqslant 6$. We have $\hat{f}(L_{i,k}^l)=\hat{f}(L_{i,k}^{l+6})=(L_{i,k}^l)^*$.

Note that $(L_{i, k}^{2l}\cup L_{i, k}^{2l+6})\cap \phi_{l, k}\neq \emptyset$, $l=1, 2$. and $(L_{i, k}^6\cup L_{i, k}^{12})\cap \phi_{3, i}\neq \emptyset$.

By Remark \ref{Lij},
\begin{equation*}
\text{The angle between $L_{j,j}^*$ and $L_{i,j}^*$, $i\neq j$, at $\hat{c}_{1,j}$ $\in$}
\begin{cases}
  (0, \pi)    & \text{$i-j\leqslant \frac{p-1}{2}$ (mod $p$)}, \\
  (\pi, 2\pi)    & \text{otherwise}.
\end{cases}
\end{equation*}
\begin{equation*}
\text{The angle between $L_{j,j}^*$ and $L_{i,j}^*$, $i\neq j$, at $\hat{c}_{2,j}$ $\in$}
\begin{cases}
  (\pi, 2\pi)   & \text{$i-j\leqslant \frac{p-1}{2}$ (mod $p$)}, \\
  (0, \pi)   & \text{otherwise}.
\end{cases}
\end{equation*}
\begin{equation*}
\text{The angle between $L_{j,j}^*$ and $L_{j,k}^*$, $j\neq k$ at $\hat{c}_{3,j}$ $\in$}
\begin{cases}
   (\pi, 2\pi)    & \text{$j-k\leqslant \frac{p-1}{2}$ (mod $p$)}, \\
    (0, \pi)      & \text{otherwise}.
\end{cases}
\end{equation*}

Now, we \textit{reorient} $\{L_{i,j}^*: 1\leqslant i,j \leqslant p\}$ by
changing the orientations of $\{L_{i,j}^*: 1\leqslant i-j\leqslant
\displaystyle{\frac{p-1}{2}}\ \text{(mod $p$)} ,1\leqslant i,j \leqslant p\}$.
We also change the orientations of the corresponding $L_{i,j}$'s.

After this reorientation, the angle between $L_{j,j}^*$ and $L_{i,j}^*$ at
$\hat{c}_{1,j}$ $\in (\pi, 2\pi)$; the angle between $L_{j,j}^*$
and $L_{i,j}^*$ at $\hat{c}_{2,j}$ $\in (0, \pi)$;  the angle
between $L_{j,j}^*$ and $L_{j,k}^*$ at $\hat{c}_{3,j}$ $\in (0, \pi)$, ($i\neq j,j\neq k$, here $i$ may equal $k$).

 If we equip $\{L_{i, k}^l: 1\leqslant i,k\leqslant p, i\neq k,
 1\leqslant l\leqslant 12\}$ with the new induced orientation, then
\begin{equation}\label{tail}
\text{the tails of} \begin{cases}
L_{i, k}^2\ \text{and}\ L_{i, k}^8\ \text{lie in}\ T_1^k, \\
L_{i, k}^4\ \text{and}\ L_{i, k}^{10}\ \text{lie in}\ T_2^k,\\
L_{i, k}^6\ \text{and}\ L_{i, k}^{12}\ \text{lie in}\ T_2^i.
\end{cases}
\text{the heads of} \begin{cases}
L_{i, k}^2\ \text{and}\ L_{i, k}^8\ \text{lie in}\ T_2^k, \\
L_{i, k}^4\ \text{and}\ L_{i, k}^{10}\ \text{lie in}\ T_1^k,\\
L_{i, k}^6\ \text{and}\ L_{i, k}^{12}\ \text{lie in}\ T_1^i.
\end{cases}
\end{equation}

The induced orientation on other segments of $L_{i, k}$'s are described in Table \ref{orientL}, where $r=0$ or $6$. We get Table \ref{orientL} from Table \ref{L*} by changing the orientation of $\{L_{i, k}^*: 1\leqslant i,k\leqslant p, 1\leqslant i- k\leqslant\frac{p-1}{2}\}$.
\bigskip

\begin{table}
\begin{center}
\newcommand{\rb}[1]{\raisebox{-1.5ex}[0pt]{#1}}
\begin{tabular}{|c|c|c||c|c|c||c|c|c|}
\hline
$\rb{$L_{i,k}^{1+r}$}$ & $\rb{tail}$ &$\rb{head}$ & $\rb{$L_{i,k}^{3+r}$}$ & $\rb{tail}$ & $\rb{head}$& $\rb{$L_{i,k}^{5+r}$}$& $\rb{tail}$ & $\rb{head}$\\[3ex]
\hline
 $\rb{${i-k}\leqslant\frac{p-1}{2}$}$ & $\rb{$T_1^i$}$ & $\rb{$T_1^k$}$ & $\rb{${i-k}\leqslant\frac{p-1}{2}$}$ & $\rb{$T_2^k$}$ & $\rb{$T_2^k$}$ & $\rb{${i-k}\leqslant\frac{p-1}{2}$}$ & $\rb{$T_1^k$}$ & $\rb{$T_2^i$}$ \\[3ex]
\hline
 $\rb{${i-k}>\frac{p-1}{2}$}$ & $\rb{$T_2^k$}$ & $\rb{$T_2^i$}$ &  $\rb{${i-k}>\frac{p-1}{2}$}$&$\rb{$T_1^k$}$ & $\rb{$T_1^k$}$ &  $\rb{${i-k}>\frac{p-1}{2}$}$&$\rb{$T_1^i$}$ & $\rb{$T_2^k$}$ \\[3ex]
\hline
\end{tabular}
\end{center}
Note: ${i-k}$ is considered as a nonnegative integer mod $p$, $r=0$ or $6$.
\caption{\label{orientL} The induced orientation on $L_{i, k}^{2l-1}, 1\leqslant l\leqslant 6$.}
\end{table}

\begin{figure}
\begin{center}
\includegraphics[width=5in]{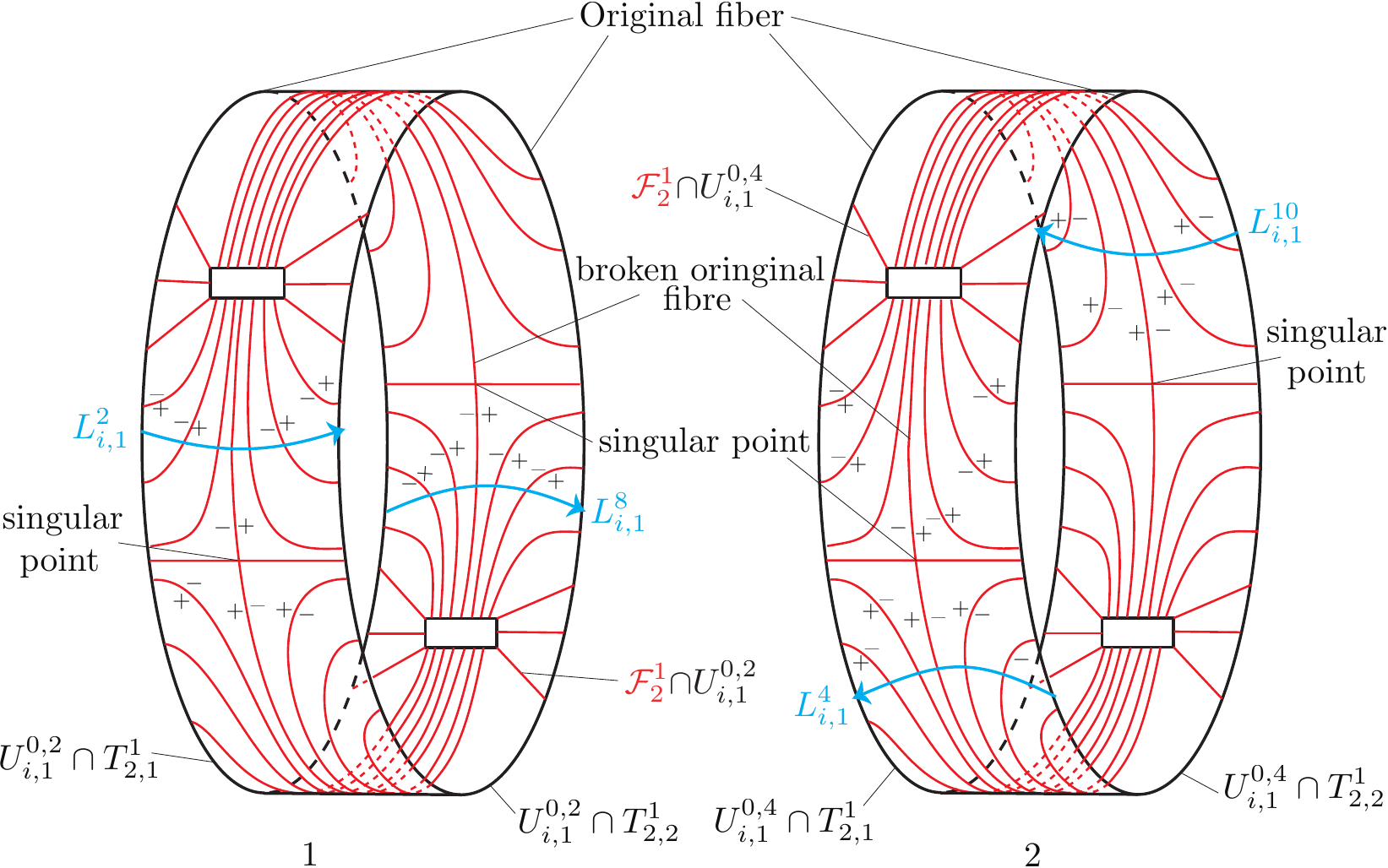}
\end{center}
\caption{\label{foliation}$U_{i, 1}^{0, 2}$ and $U_{i, 1}^{0, 4}$.}
\end{figure}

Fix a transverse orientation on the surface bundle $\mathcal{F}_2^j$ in $M_2^j$.
By the construction of $H$ in Step 2, there are two
 singular points in the foliation of $U_{i,k}^{0,2l}$ given by
 $U_{i,k}^{0,2l}\cap\mathcal{F}_2$ if $\delta$ is small enough, $l=1,2,3$ (cf. Figure \ref{foliation}).
 We can isotope $L_{i,k}^{2l+r}$ in $U_{i,k}^{0,2l}$
such that $L_{i,k}^{2l}$ is transverse to
$\overset{p}{\underset{j=1}{\cup}}\mathcal{F}_2^j$, and travels
from the ``$-$" side of $\mathcal{F}_2^j$ to the ``$+$" side
($1\leqslant l\leqslant3,\ r=0,6$). We only discuss the isotopy in $M_2^1$,
and others are similar.

By the definition of $U_{i,k}^{0,l}(1\leqslant l\leqslant 6)$, ($U_{i,1}^{0,2}\cup U_{i,1}^{0,4}\cup U_{1,i}^{0,6})\subset M_2^1$, $2\leqslant i\leqslant p$. Suppose that $\delta$ is small enough so that  the singular foliation on $U_{i, 1}^{0, 2}$, $U_{i, 1}^{0, 4}$ or $U_{1, i}^{0, 6}$ induced by $\mathcal{F}_2^1$ has two singular points. 
The singular points on $U_{i,1}^{0,2l} \in \phi_{l, 1}$ $(l=1,2)$, and the singular points on $U_{1,i}^{0,6} \in \phi_{3, 1}$, $2\leqslant i\leqslant p$.
Recall that $L_{i,k}^l\cup L_{i,k}^{l+6}\in U_{i,k}^{0,l}$, $1\leqslant l \leqslant 6$,  so $\{L_{i,k}^{2l+r}: 1\leqslant i,k\leqslant p, i\neq k, (k,l)=(1,1), (1,2), \text{ or } (i,l)=(1,3), r=0,6\} =L\cap M_2^1$. 

 Let $S=\{(i,k,l,r): 1\leqslant i,k \leqslant p, i\neq k, (k,l)=(1,1), (1,2)$, or $(i,l)=(1,3), r=0,6\}$. Fix $(i,k,l,r)\in S$. By Lemma 2.1 in \cite{abz}, $L_{i,k}^{2l+r}$ is always perpendicular to the original fibers in $U_{i,k}^{0,2l}$.
If $L_{i,k}^{2l+r}$ passes from the ``$+$'' to the ``$-$'' side of $\mathcal{F}_2^1$'s leaves along its orientation, we isotope $L_{i,k}^{2l+r}$ passing the singular point smoothly along the original fibres  in a small neighborhood of $U_{i,k}\cap M_2^1$, and place them as in Figure \ref{foliation}.
By (\ref{tail}), $L_{i, 1}^{2+r}\subset U_{i, 1}^{0, 2}$ travels from $T_{2, 1}^1$ to $T_{2, 2}^1$; $L_{i, 1}^{4+r}\subset U_{i, 1}^{0, 4}$ travel from $T_{2, 2}^1$ to $T_{2, 1}^1$. Then $L_{i, 1}^2$ and $L_{i, 1}^4$ are in the different sides of their own singular points, so are $L_{i, 1}^8$ and $L_{i, 1}^{10}$. 
Also by (\ref{tail}), $L_{i, 1}^{4+r}\subset U_{i, 1}^{0, 4}$ and $L_{1, i}^{6+r}\subset U_{1, i}^{0, 6}$ are both travel from $T_{2, 2}^1$ to $T_{2, 1}^1$, so $U_{1, i}^{0, 6}\cap \mathcal{F}_2^1$ is similar to $U_{i, 1}^{0, 4}\cap \mathcal{F}_2^1$. $U_{i, 1}^{0, 2}\cap \mathcal{F}_2^1$ is as shown in Figure \ref{foliation}-1, and $U_{i, 1}^{0, 4}\cap \mathcal{F}_2^1$ is as shown in Figure \ref{foliation}-2, in other words, they alternate with the singular points in the original circle fiber direction.

\begin{figure}
\begin{center}
\includegraphics[width=5.9in]{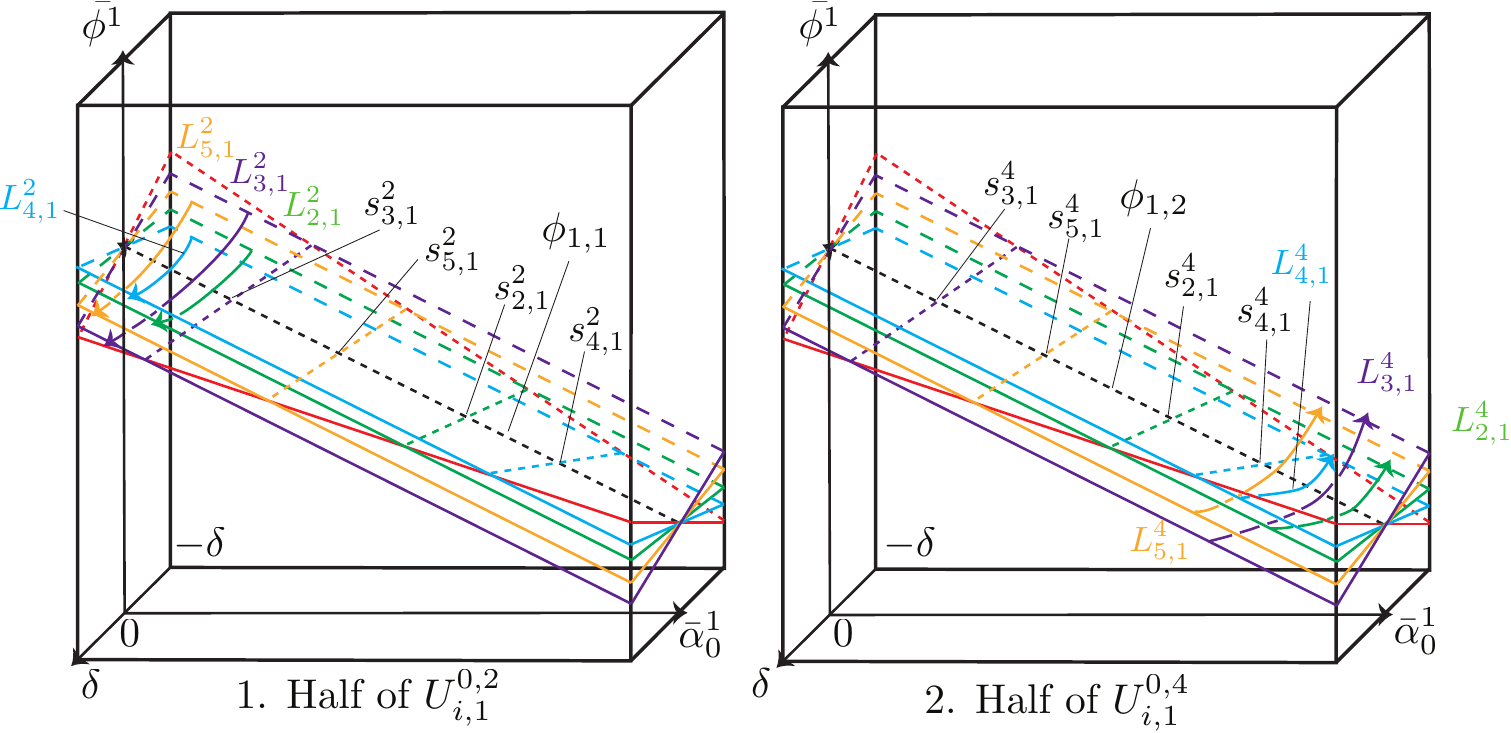}
\end{center}
\caption{\label{isotope}The positions of $L_{i, 1}^2$ and $L_{i, 1}^4$.}
\end{figure}

The isotopies are similar to that in \cite{abz}. Next, we want to show that the above isotopies of $L_{i,k}^{2l+r}$ won't block each other. Since $U_{i,1}^{0,2}\cap U_{k,1}^{0,2}=\phi_{1,1}$, $U_{i,1}^{0,4}\cap U_{k,1}^{0,4}=\phi_{2,1}$, and $U_{1,i}^{0,6}\cap U_{1,k}^{0,6}=\phi_{3,1}$, $2\leqslant i,k \leqslant p, i\neq k$, potential obstructions happen  in the set $\{L_{i,1}^{2+r}: 2\leqslant i\leqslant p, r=0,6\}$, $\{L_{i,1}^{4+r}: 2\leqslant i\leqslant p, r=0,6\}$, or $\{L_{1,i}^{6+r}: 2\leqslant i\leqslant p, r=0,6\}$.  
$L_{i,k}^{2l}$ and $L_{i,k}^{2l+r}$ are in the different halves of $U_{i,k}^{0,2l}$ (cf. Figure \ref{foliation}), so they won't block each other since the smooth $\phi$-vertical isotopy can not cross the punctures on $U_{i,k}^{0,2l}$.  
Next we   check $\{L_{i,1}^{2}: 2\leqslant i\leqslant p\}$, and the discussions for $\{L_{i,1}^{4}: 2\leqslant i\leqslant p\}$ and $\{L_{1,i}^{6}: 2\leqslant i\leqslant p\}$ are similar. 
By (\ref{tail}), $L_{i, 1}^{2}$ travels from $T_{2, 1}^1$ to $T_{2, 2}^1$
for any $2\leqslant i \leqslant p$. We may assume that we need isotope $L_{i,1}^2$ smoothly along $\phi$ above its singular point, as shown in Figure \ref{isotope}-1, $2\leqslant i \leqslant p$. (There is only half of $U_{i,1}^{0,2}\cap (T^1\times [-\delta,\delta])$ in Figure \ref{isotope}-1). Although the singular points of $U_{i,1}^{0,2}$'s in $\phi_{1,1}$ are at different positions, and the distribution of the intersection points of  $\phi_{1,1}$ and $\{L_{i,1}^2:2\leqslant i \leqslant p\}$ may be different for different $(q_1, q_2,q_3)$, we can isotope $L_{i,1}^2$ in the positive direction of $\phi_{1,1}$ above all singular points, so they won't block each other. Figure \ref{isotope}-2 illustrates the situation near $U_{i,1}^{0,4}\cap T^1\times[-\delta,\delta]$. In Figure \ref{isotope}, $s_{i, 1}^{2l}$ denotes the singular point of $L_{i, 1}^{2l}$, for $i=2, 3, 4, 5$ and $l=1, 2$.

 By generalizing  the above discussion  to $M_2^j$, $1\leqslant j \leqslant p$, we have the following proposition which is analogous to Proposition 6.2 in \cite{abz}.

\begin{prop}\label{M2j}
Suppose $1\leqslant j \leqslant p$. Fix a transverse orientation on $\mathcal{F}_2^j$ in $M_2^j$. There is a smooth $\phi$-vertical isotope of $L_{i, j}$ and $L_{j, i}$ $(1 \leqslant i\leqslant p, i\neq j)$ in $\overset{p}{\underset{i=1,i\neq j}{\cup}}(U_{i, j}^0\cup U_{j, i}^0)$, fixed outside a small neighborhood of $(\overset{p}{\underset{i=1,i\neq j}{\cup}}(U_{i, j}^0\cup U_{j, i}^0))\cap M_2^j)$, which reposition $L_{i, j}^l$ and $L_{j, i}^k$ to be transverse to $\mathcal{F}_2^j$ and to pass from the negative to the positive side of $\mathcal{F}_2^j$'s leaves while traveling along the orientation on $L_{i, j}^l$ and $L_{1, j}^k$, where $1\leqslant i\leqslant p, i\neq j; l=2, 8, 4, 10; k=6, 12$. 
\end{prop}

Proposition \ref{M2j} finishes Step 3.

\noindent\textbf{Step 4}

\begin{figure}
\begin{center}
\includegraphics[width=4in]{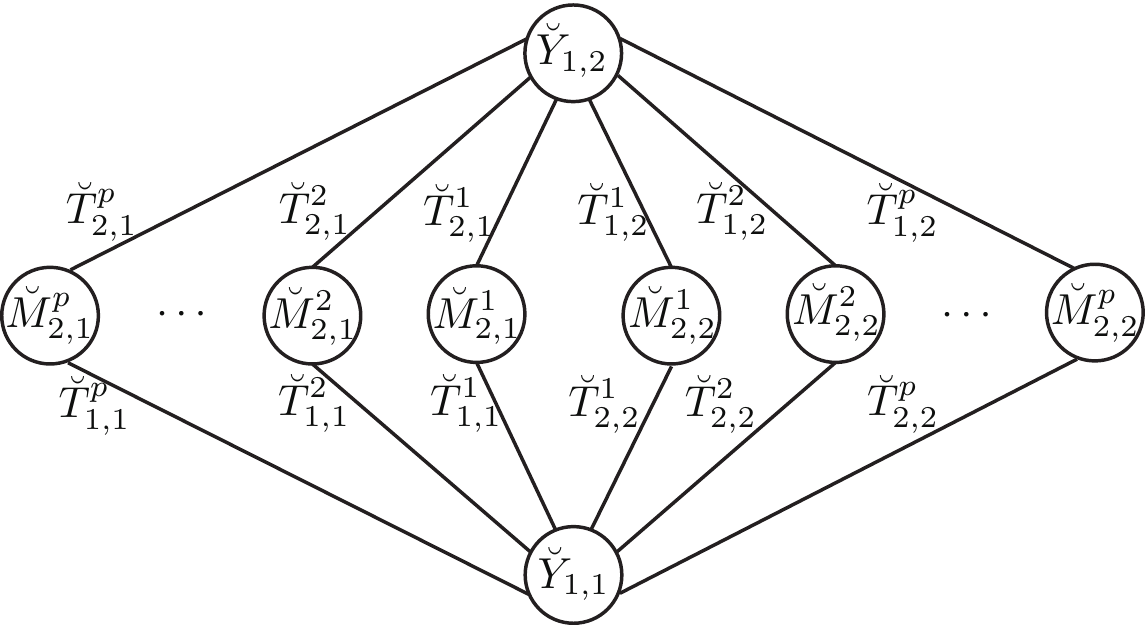}
\end{center}
\caption{\label{double}The graph decomposition of $\breve{M}$.}
\end{figure}

At first we construct a double cover of $M$ similarly to \cite{abz}.  Set $M_0$ to be the 3-manifold obtained by cutting $M$ open along $\{T_2^j$, $1\leqslant j \leqslant p\}$. 
Take two copies of $M_0$, denoted $\breve{M}_{0,1}$ and $\breve{M}_{0,2}$. We have $\breve{M}_{0,s}=\breve{Y}_{1,s}\cup \breve{M}_{2,s}^1\cup \breve{M}_{2,s}^2\cup \cdots\cup \breve{M}_{2,s}^p$, where $\breve{Y}_{1,s}$ is a copy of $Y_1$, and $\breve{M}_{2,s}^j$ is a copy of $M_2^j$, $1\leqslant j\leqslant p$ and $s=1,2$. 
$\breve{M}$ is a 3-manifold obtained by gluing $\breve{M}_{2,1}^j$ to $\breve{Y}_{1,2}$ along $\breve{T}_{2,1}^j$, and $\breve{M}_{2,2}^j$ to $\breve{Y}_{1,1}$ along  $\breve{T}_{2,2}^j$, where $\breve{T}_{2,s}^j\subset \breve{M}_{0,s}$  is a copy of $T_2^j$, $1\leqslant j\leqslant p$ and $s=1,2$. 
The gluing map is the same as the one used to glue $M_2^j$ back to $Y_1$ to get $M$ for every $1\leqslant j\leqslant p$. This construction is shown in Figure \ref{double}.
It's clear that $\breve{M}$ is a free double cover of $M$. We denote $p_2$ to be the covering map. Then
\begin{equation*}
p_2^{-1}(Y_1)=\breve{Y}_{1, 1}\cup \breve{Y}_{1, 2}
\end{equation*}
\begin{equation*}
p_2^{-1}(M_2^j)=\breve{M}_{2, 1}^j\cup \breve{M}_{2, 2}^j
\end{equation*}
\begin{equation*}
p_2^{-1}(T_r^j)=\breve{T}_{r, 1}^j\cup \breve{T}_{r, 2}^j\ (r=1, 2)
\end{equation*}
 $\breve{M}$ is also a graph manifold and Figure \ref{double}  illustrates the graph of the JSJ-decomposition of $\breve{M}$.

By construction, $\breve{M}$ is fibred. In Step 2, we construct $H=H_1\cup(\overset{p}{\underset{j=1}{\cup}}H_2^j)$ as a leaf of the semi-bundle structure of $M$. 
Let $\breve{H}=p_2^{-1}(H)$. It's not hard to see that $\breve{H}$ is a connected orientable surface in $\breve{M}$. 
Let 
\begin{equation*}
\breve{\mathcal{F}}=\breve{\mathcal{F}}_{1, 1}\cup \breve{\mathcal{F}}_{1, 2}\cup (\underset{s=1}{\overset{2}{\cup}}\underset{j=1}{\overset{p}{\cup}}\breve{\mathcal{F}}_{2, s}^j),
\end{equation*}
 where $\breve{\mathcal{F}}_{1, s}=p_2^{-1}(\mathcal{F}_1)\cap \breve{Y}_{1, s}, \breve{\mathcal{F}}_{2, s}^j=p_2^{-1}(\mathcal{F}_2^j)\cap \breve{M}_{2, s}^j, s=1, 2$.
Then $\breve{F}$ is a surface bundle in $\breve{M}$, and $\breve{H}$ is one leaf of $\breve{F}$.  Fix a transverse orientation for $\breve{\mathcal{F}}$ and let $\breve{\mathcal{F}}_{1, s}$ and $\breve{\mathcal{F}}_{2, s}^j$ have the induced transverse orientation, $s=1, 2$, $1\leqslant j\leqslant p$.

$p_2$ can be extended to a double cover of $Y$. Denote the covering space by $\breve{Y}$. $\breve{Y}$ inherits the Seifert fibred structure from $Y$ with fibre $\breve{\phi}$. Let $\breve{\phi}$ have the inherited orientation from $\phi$.
\begin{equation*}
\breve{Y}=\breve{Y}_{1, 1}\cup \breve{Y}_{1, 2}\cup (\underset{s=1}{\overset{2}{\cup}}\underset{j=1}{\overset{p}{\cup}}\breve{Y}_{2, s}^j)
\end{equation*}
 where $\breve{Y}^j_{2, 1}\cup\breve{Y}_{2,2}^j=p_2^{-1}(Y_2^j)$ and $\breve{M}_{2, s}^j\subset\breve{Y}_{2,s}^j$, $s=1, 2$, $1\leqslant j\leqslant p$. By construction, $p_2^{-1}(L_{j,j})$ are two copies of $L_{j,j}$, denoted $\breve{L}_{j,j,1}$, and $\breve{L}_{j,j,2}$. $\breve{L}_{j,j,s}\subset{\breve{Y}_{2,s}^j}$, $s=1,2$. Then 
 \begin{equation*}
 \breve{M}=\breve{Y}\setminus (\underset{s=1}{\overset{2}{\cup}}\underset{j=1}{\overset{p}{\cup}} N(\breve{L}_{j,j,s})).
\end{equation*}
Equip $\breve{L}_{j,j,s}$ with the inherited orientation from $L_{j,j}$, $1\leqslant j\leqslant p, s=1,2$.

By the construction of $\breve{M}$, $p_2^{-1}(U_{i,k})$ is a $\breve{\phi}$-vertical torus, $1\leqslant i,k\leqslant p, i\neq k$. Let $\breve{U}_{i, k}=p_2^{-1}(U_{i, k})$, and $\breve{U}_{i,k}^0=\breve{U}_{i,k}\cap\breve{M}$. 
Denote $p_2^{-1}(U_{i, k}^{0, l})=\breve{U}_{i, k, 1}^{0, l}\cup \breve{U}_{i, k, 2}^{0, l}$, $1\leqslant l\leqslant 6$, where $(\breve{U}_{i, k, s}^{0, 1}\cup \breve{U}_{i, k, s}^{0, 3}\cup \breve{U}_{i, k, s}^{0, 5})\subset \breve{Y}_{1, s}$, $(\breve{U}_{i, k, s}^{0, 2}\cup \breve{U}_{i, k, s}^{0, 4})\subset \breve{M}_{2, s}^k$, and $\breve{U}_{i, k, s}^{0, 6}\subset \breve{M}_{2, s}^i, s=1, 2$.
Denote $\breve{L}=p_2^{-1}(L)$. $p_2^{-1}(L_{i,k})$ has two components, and two-fold trivially covers its Seifert quotient image, since $L_{i,k}$ is a two-fold cyclic cover of $L^*_{i,k}$. Then $\breve{L}$ has $2p^2$ components. 
Suppose $(\breve{L}_{i, k, s}^l\cup \breve{L}_{i, k, s}^{l+6})\subset \breve{U}_{i, k, s}^{0, l}$ to be the lift of $(L_{i, k}^l\cup L_{i, k}^{l+6})\subset U_{i, k}^{0, l}, 1\leqslant l\leqslant 6$. Equip $\breve{L}_{i, k, s}^l$ with the inherited orientation, $1\leqslant l\leqslant 12, s=1, 2$.

Now we fix $s=1, 2$ in the following discussion.

By construction, proposition \ref{M2j} holds for every $\breve{M}_{2, s}^j$, $1\leqslant j\leqslant p$. In particular, we can isotope $\{\breve{L}_{i, j,s}\cup \breve{L}_{j, i,s}: 1 \leqslant i,j\leqslant p, i\neq j, s=1,2\}$ smoothly along $\breve{\phi}$ such that $\{\breve{L}^{2+r}_{i, j,s}\cup\breve{L}^{4+r}_{i,j,s}\cup\breve{L}^{6+r}_{j, i,s}: 1 \leqslant i,j\leqslant p, i\neq j, r=0,6, s=1,2\}$ pass from the negative side to the positive side of $\breve{\mathcal{F}}_{2, s}^j$'s leaves when traveling along their orientation  in $\breve{M}_{2, s}^j$. This isotopy is fixed outside a small neighborhood of $\overset{2}{\underset{s=1}{\cup}}(\overset{p}{\underset{i=1,i\neq j}{\cup}}(U_{i, j,s}^0\cup U_{j, i,s}^0))\cap\breve{M}_{2,s}^j)$.
Then $\{\breve{L}_{i,k,s}: 1\leqslant i,k\leqslant p, i\neq k, s=1,2\}$ are transverse to $\{\breve{\mathcal{F}}_{2,s}^j, 1\leqslant j \leqslant p, s=1,2\}$.  
Next, we perform Dehn twist operations a sufficiently large number of times on $\{\breve{\mathcal{F}}_{1, s} : s=1, 2\}$ along a set of tori $\Gamma$, such that the new surface fiber is transverse to $\{\breve{L}_{i, k, s}, 1\leqslant i,k\leqslant p, i\neq k, s=1,2\}$ everywhere. Unlike in \cite{abz}, here we need more vertical tori. $\Gamma$ consists $2p+4$ $\breve{\phi}$-vertical tori,$\{\breve{V}_s^j: 1\leqslant j \leqslant p+2, s=1,2\}$. $\{\breve{V}_s^j: 1\leqslant j \leqslant p, s=1,2\}$ are $2p$ vertical tori, where $\breve{V}_s^j$ parallel to $\breve{T}_{1, s}^j$ for $1\leqslant j \leqslant p, s=1,2$. In addition to these boundary parallel tori, we need $4$ extra tori, which we will give the construction later.

\begin{figure}
\begin{center}
\includegraphics[width=6in]{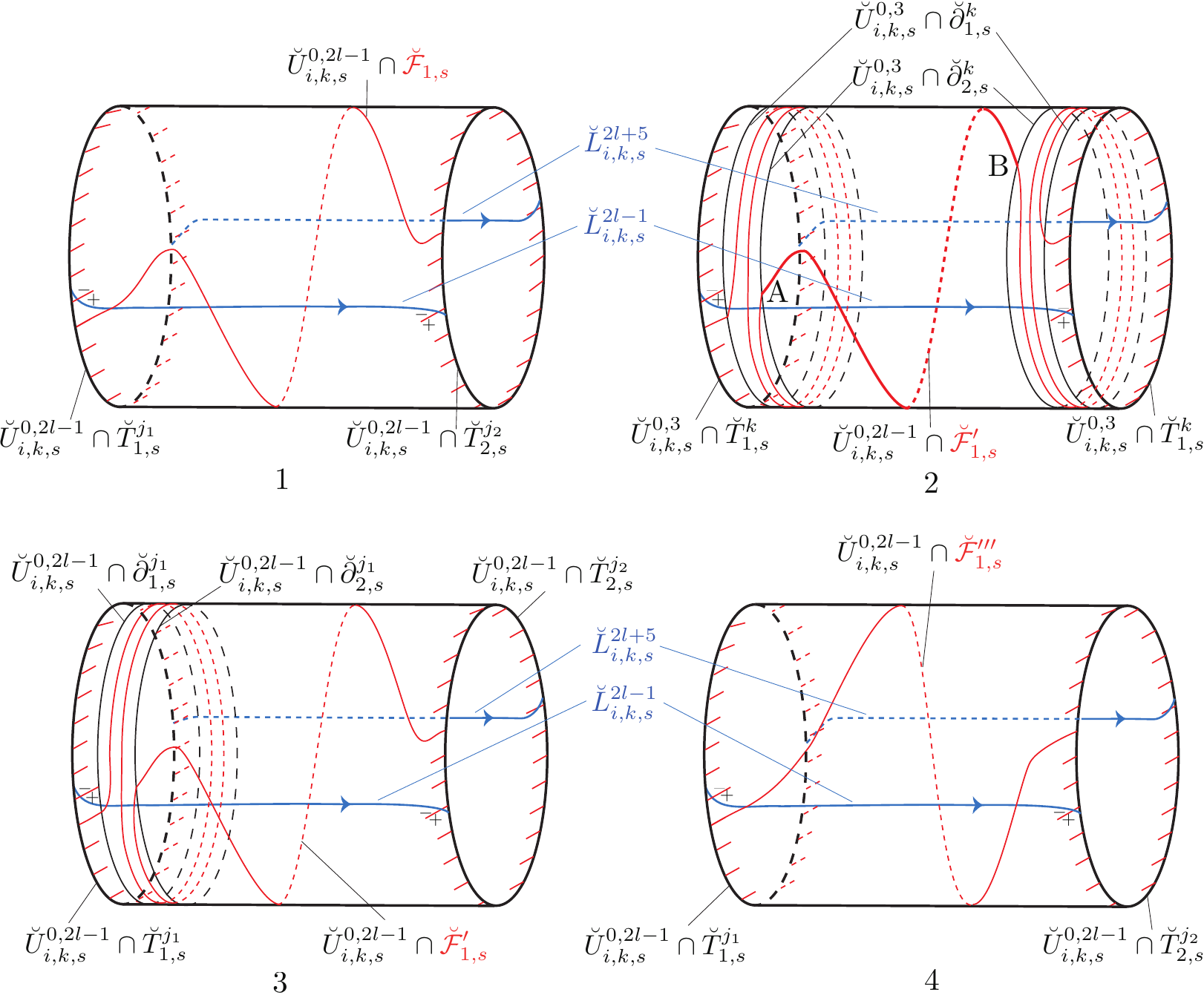}
\end{center}
\caption{\label{dehntwistU}Dehn twist of $\breve{\mathcal{F}}_{1, s}$.}
\end{figure}

$\breve{U}_{i,k,s}^{0,2l-1}$ is a $\breve{\phi}$-vertical annulus, $1\leqslant l\leqslant 6$. $\breve{L}_{i,k,s}^{2l-1}\cup\breve{L}_{i,k,s}^{2l+5}\subset\breve{U}_{i,k,s}^{0,2l-1}$, and transverse to $\breve{\phi}$ by Lemma 2.1 in \cite{abz}, $1\leqslant l\leqslant 3$. By the construction of $\mathcal{F}$ in Step 2, we may assume that the oriented arcs $\breve{L}_{i,k,s}^{2l-1}\cup\breve{L}_{i,k,s}^{2l+5}$ are always transverse from the negative to the positive side of $\breve{\mathcal{F}}_{1,s}$ near the boundary of $\breve{U}_{i,k,s}^{0,2l-1}$, $1\leqslant l\leqslant 3$. By the construction of $H_1$ in Step 2, $\breve{\mathcal{F}}_{1,s}$ is transverse to the original fibre $\breve{\phi}$. $\breve{\mathcal{F}}_{1,s}\cap\breve{U}_{i,k,s}^{0,2l-1}$ gives another foliation of $\breve{U}_{i,k,s}^{0,2l-1}$, $1\leqslant l\leqslant 3$. One leaf of $\breve{\mathcal{F}}_{1,s}\cap\breve{U}_{i,k,s}^{0,2l-1}$ is shown in Figure \ref{dehntwistU}.

 Let $\breve{V}_s^j$ be a $\breve{\phi}$-vertical torus in a small neighborhood of $\breve{T}_{1, s}^j$ in $int(\breve{Y}_{1, s})$ and parallel to $\breve{T}_{1, s}^j$, $1\leqslant j\leqslant p, s=1,2$. By the above discussion we may assume that $\breve{V}_s^j$ intersect $\breve{L}_{i,k,s}^{2l-1}\cup\breve{L}_{i,k,s}^{2l+5}$ transversely, $1\leqslant l\leqslant 3$. 
 
\begin{figure}
\begin{center}
\includegraphics{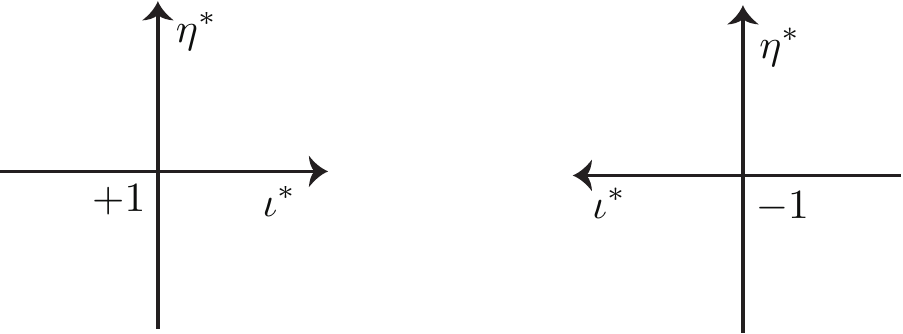}
\end{center}
\caption{\label{sign}Algebraic intersection number}
\end{figure}

 Now we introduce a sign for every intersection point of an arc $\eta$ and a $\breve{\phi}$-vertical torus $T$ in $\breve{M}$. Suppose $\eta^*$ and $\iota^*$ are the  images of $\eta$ and $T$ under the Seifert quotient map respectively. We say $\eta$ intersects T \textit{positively} if the algebraic intersection number $i(\eta^*, \iota^*)= +1$(cf. Figure \ref{sign}), otherwise we say $\eta$ intersects T \textit{negatively}. Equip the image of $\breve{V}_s^j$ under the Seifert quotient map the same orientation as $\breve{L}_{j,j,s}$, $1\leqslant j\leqslant p$. It's easy to see that if the tail of $\breve{L}_{i,k,s}^{2l-1+r}$ is in $\breve{T}_{1,s}^j$, $\breve{L}_{i,k,s}^{2l-1+r}$ intersects $\breve{V}_s^j$ positively,(cf. Figure \ref{dehntwistU}-1),$1\leqslant l\leqslant 3, r=0,6$.

The following remark is obtained from Table \ref{orientL}. 

\begin{remark}\label{intersection}
Any of $\{\breve{L}_{i,k, s}^{1+r}, {i-k} \leqslant\displaystyle\frac{p-1}{2},r=0,6,s=1,2\}\cup\{\breve{L}_{i,k,s}^{3+r}, {i-k}>\displaystyle\frac{p-1}{2}, r=0,6,s=1,2\}$, intersects $\overset{2}{\underset{s=1}{\cup}}\overset{p}{\underset{j=1}{\cup}}\breve{V}_s^j$ twice and in different signs, and any of $\{\breve{L}_{i, k, s}^{5+r} : 1\leqslant i,k \leqslant p, i\neq k, r=0,6, s=1,2\}$ intersects $\overset{2}{\underset{s=1}{\cup}}\overset{p}{\underset{j=1}{\cup}}\breve{V}_s^j$ once and positively. Note that $\overset{2}{\underset{s=1}{\cup}}\overset{p}{\underset{j=1}{\cup}}\breve{V}_s^j$ only intersects the above arcs.
\end{remark}

By the above remark, when we perform Dehn twist along $\{\breve{V}_s^j: 1\leqslant j \leqslant p, s=1,2\}$, the relation between $\{\breve{L}_{i,k, s}^{1+r}, {i-k} \leqslant\displaystyle\frac{p-1}{2},r=0,6,s=1,2\}\cup\{\breve{L}_{i,k,s}^{3+r}, {i-k}>\displaystyle\frac{p-1}{2}, r=0,6,s=1,2\}$ and the new fibres are the same as the relation of those with the old fibres. To make $\{\breve{L}_{i,k, s}^{1+r}, {i-k} \leqslant\displaystyle\frac{p-1}{2},r=0,6,s=1,2\}\cup\{\breve{L}_{i,k,s}^{3+r}, {i-k}>\displaystyle\frac{p-1}{2}, r=0,6,s=1,2\}$ transverse to the new fibres, we need extra Dehn twists. We now construct $\{\breve{V}_s^{p+1}\cup \breve{V}_s^{p+2}\subset \breve{Y}_{1,s},s=1,2\}$. $\breve{V}_s^{p+r}$ is the lift of $\phi$-vertical torus $V^{p+r} \subset Y_1$, where  $V^{p+r}=\hat{f}^{-1}(l_r), r=1,2$. $l_1$ and $l_2$ are two simple closed curves in $F_1$ and constructed as follows.

\begin{figure}
\begin{center}
\includegraphics[width=5.9in]{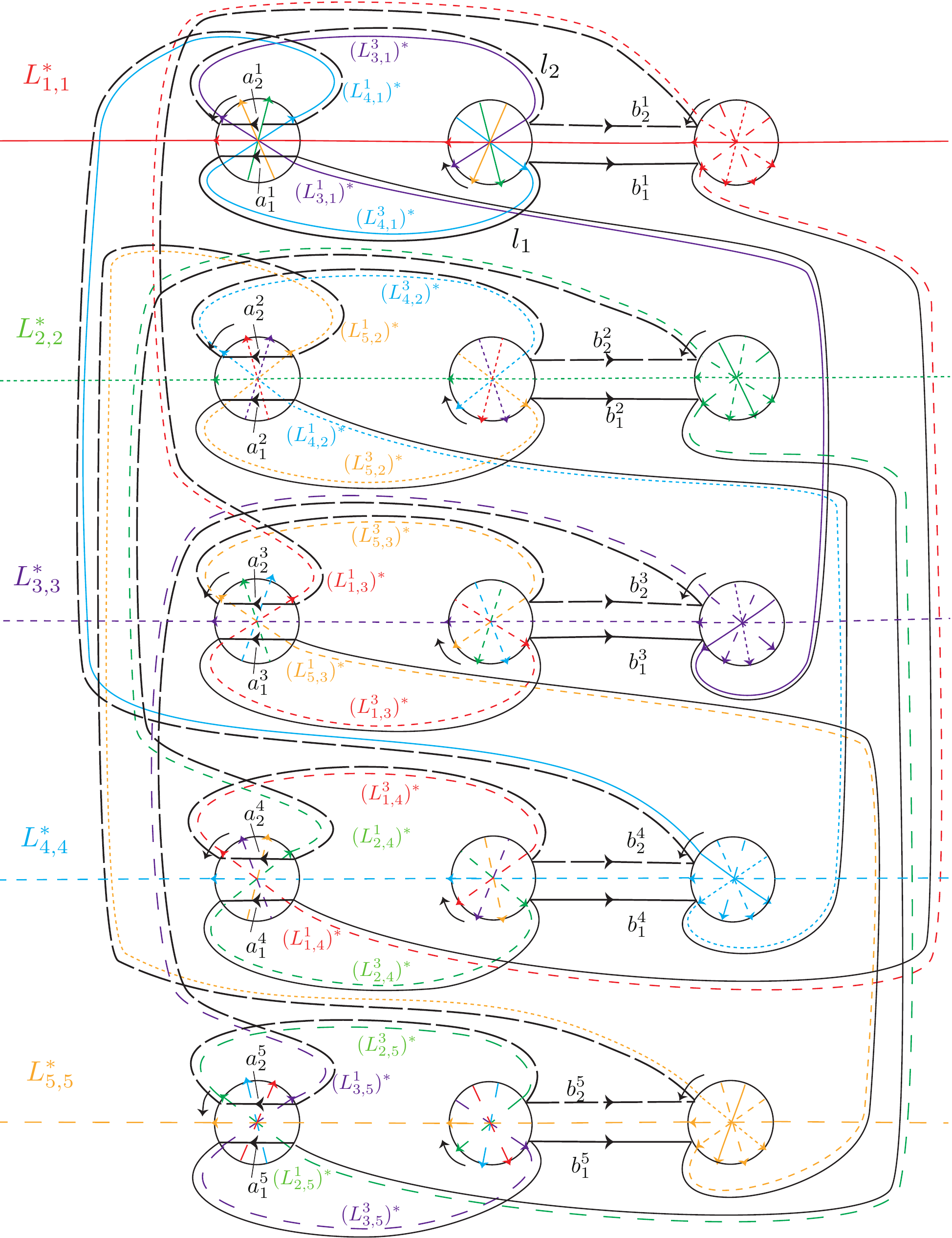}
\end{center}
\caption{\label{l1l2}$l_1$ and $l_2$}
\end{figure}

Recall that $\beta_1^j=L_{j,j}^*\times\{-\epsilon\}$, $\beta_2^j=L_{j,j}^*\times\{\epsilon\}$, $1\leqslant j\leqslant p$.
Let $a_r^j$ be the intersection of $\beta_r^j$ and a small disk centered at $\hat{c}_{1,j}$ with radius bigger than $\epsilon$, $1\leqslant j\leqslant p$. Push $a_r^j$ into $F_1$ a little bit, as shown in Figure \ref{l1l2}. By Remark \ref{Lij}, $\{L_{i,j}^*, 1\leqslant i, j\leqslant p, i\neq j\}$ intersect $L_{j,j}^*$ at $\hat{c}_{1,j}$. Then $a_r^j\cap L_{i,j}^*\neq\emptyset$, $1\leqslant i,j\leqslant p$, $i\neq j$, $r=1,2$.
Since $\{L_{i,k}^* : 1\leqslant i, k\leqslant p, i\neq k\}$ intersects $L_{j,j}^*$ only at $\hat{c}_{1,j}$,  $\hat{c}_{2,j}$ and $\hat{c}_{3,j}$, we can take an arc $b_r^j$ in $F_1$ which is parallel to $\beta_r^j$ and between $\hat{c}_{2,j}$ and $\hat{c}_{3,j}$, $1\leqslant j\leqslant p$, $r=1,2$ (cf. Figure \ref{l1l2}).

By Remark \ref{Lij} and Remark \ref{arcs}, we can connect
\begin{equation}\label{connect}
\begin{cases}
   a_1^j \ \text{and} \ b_1^j,   & \text{ by the arc parallel to $(L_{j+\frac{p+1}{2},j}^3)^*$}, \\
    b_1^j \ \text{and} \ a_1^{j+\frac{p+1}{2}}  & \text{by the arc parallel to $(L_{j,j+\frac{p+1}{2}}^1)^*$};\\
    a_2^j \ \text{and} \ b_2^j,   & \text{ by the arc parallel to $(L_{j+\frac{p-1}{2},j}^3)^*$}, \\
    b_2^j \ \text{and} \ a_2^{j+\frac{p-1}{2}}  & \text{by the arc parallel to $(L_{j,j+\frac{p-1}{2}}^1)^*$}.\\
\end{cases}
\end{equation}
$1\leqslant  j\leqslant p$. Note that the indices are in $\mathbb{Z}/p$.
Then we get two simple closed curve $l_1$, and $l_2$ in $F_1$, (cf. Figure \ref{l1l2}).
$l_1$ and $l_2$ can be oriented as in Figure \ref{l1l2}, such that the orientations on $a_1^j$ and $a_2^j$ are the same as the orientation on $L_{j,j}^*$.

By the construction of $l_1$ and $l_2$, they only intersect $L_{i,k}^*$ at the arcs $a_1^j$ and $a_2^j$, $1\leqslant j\leqslant p$. $\{a_1^j\cup a_2^j: 1\leqslant j \leqslant p\}$ only intersects the arcs of $L_{i,k}^*$ between $\hat{c}_{3,i}$ and $\hat{c}_{1,k}$, and between $\hat{c}_{1,k}$ and $\hat{c}_{2,k}$, i.e. $(L_{i,k}^1)^*$ and $(L_{i,k}^3)^*$ by Remark \ref{arcs}. $l_1\cup l_2 $ intersects $(L_{i,k}^1)^*$ or $(L_{i,k}^3)^*$ exactly once and does not intersect $(L_{i,k}^5)^*$, $1\leqslant i,k\leqslant p$, $i\neq k$. Moreover $i((L_{i,k}^1)^*, l_r)=i((L_{i,k}^3)^*,l_r)=-1, r=1,2$, since the angle between $L_{j,j}^*$ and $L_{i,j}^*$ at $\hat{c}_{1,j}\in (\pi,2\pi)$ after the reorientation in Step 3(cf. Figure \ref{l1l2}). 
Since $\breve{V}_s^{p+m}$ and $\breve{L}_{i,k,s}^{2m-1+r}$ are the lifts of $V^{p+m}$ and $L_{i,k}^{2m-1+r}$ respectively, and their images under the Seifert quotient map are $l_m$ and $(L_{i,k}^{2m-1})^*$, respectively, $m=1,2,r=0,6$, we have the following remark.

\begin{remark}\label{intersection2}
$\overset{2}{\underset{s=1}{\cup}}(\breve{V}_s^{p+1}\cup \breve{V}_s^{p+2})$ intersects $\breve{L}_{i,k,s}^{1+r}$ and $\breve{L}_{i,k,s}^{3+r}$ exactly once and negatively, and  $\overset{2}{\underset{s=1}{\cup}}(\breve{V}_s^{p+1}\cup \breve{V}_s^{p+2})\cap\breve{L}_{i,k,s}^{5+r}=\emptyset$, $1\leqslant i,k\leqslant p$, $i\neq k$, $r=0,6.$
\end{remark}

Note that $l_1\cup l_2 \cap \{L_{j,j}^*, 1\leqslant j \leqslant p\}=\emptyset$, so we may suppose that the tori in $\Gamma$ are mutually disjoint.

Let $N(\breve{V}_s^j)$ be a small regular neighborhood of $\breve{V}_s^j$ in $int(\breve{Y}_{1, s})$ consisting of $\breve{\phi}$-circle fibres, $1\leqslant j\leqslant p+2$. We may suppose that $N(\breve{V}_s^i)\cap N(\breve{V}_s^j)=\emptyset\ (1\leqslant i, j\leqslant p+2; i\neq j)$ if we take the regular neighborhoods small enough. Define $\partial N(\breve{V}_s^j)=\breve{\partial}_{1, s}^j\cup \breve{\partial}_{2, s}^j, 1\leqslant j\leqslant p+2$. $\breve{\partial}_{1, s}^j$ and $\breve{\partial}_{2, s}^j$ are two $\breve{\phi}$-vertical tori, $1\leqslant j\leqslant p+2$.

By Remark \ref{intersection} $\overset{p}{\underset{j=1}{\cup}}N(\breve{V}_s^j)$ intersects each of $\breve{L}_{i, k, s}^5$ and $\breve{L}_{i, k, s}^{11}$ in a single arc, $1\leqslant i,k \leqslant p$, $i\neq k$, $s=1,2$. We may assume that the tail of this arc (with the induced orientation) is contained in $\breve{\partial}_{1, s}^j$, since $\{\breve{L}_{i,k,s}^5\cup\breve{L}_{i,k,s}^{11}:1\leqslant i,k\leqslant p,i\neq k, s=1,2\}$ intersects $\overset{p}{\underset{j=1}{\cup}}N(\breve{V}_s^j)$ in the same sign.
On the other hand, $N(\breve{V}_s^{p+1})\cup N(\breve{V}_s^{p+2})$ intersects each of $\breve{L}_{i, k, s}^1, \breve{L}_{i, k, s}^7, \breve{L}_{i, k, s}^3, \breve{L}_{i, k, s}^9$ in a single arc, by Remark \ref{intersection2}. We also assume that the tail of this arc (with the induced orientation) contained in $\breve{\partial}_{1, s}^{p+r}$, $r=1,2$.

\begin{figure}
\begin{center}
\includegraphics[width=4.3in]{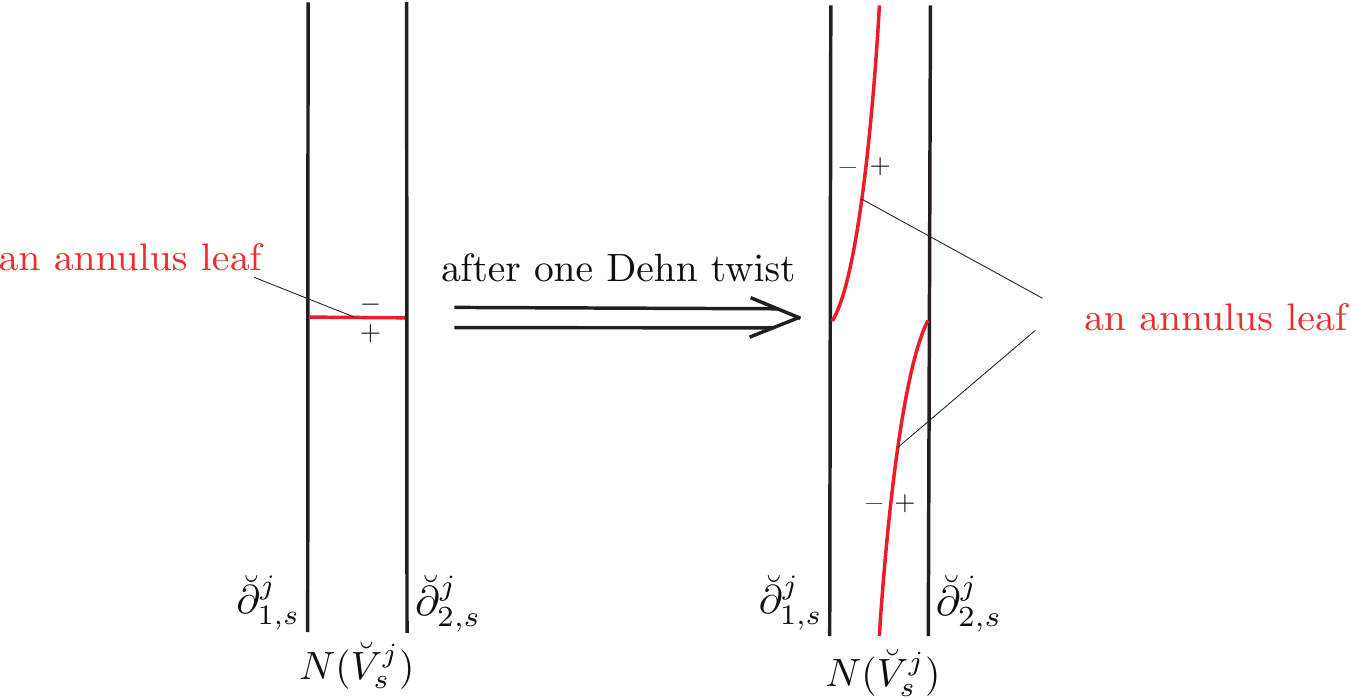}
\end{center}
\caption{\label{DT}Dehn twist in $N(\breve{V}_s^j)$ in the direction opposite to  the transverse orientation of $\breve{\mathcal{F}}_{1, s}$.}
\end{figure}

Since $\breve{\mathcal{F}}_{1, s}$ is transverse to the old fibres, $\breve{\mathcal{F}}_{1, s}\cap N(\breve{V}_s^j)$ is a foliation of $N(\breve{V}_s^j)$ by annuli, $1\leqslant j\leqslant p+2$, $s=1,2$.
The Dehn twist operation on  $\breve{\mathcal{F}}_{1, s}$, denoted by $D_s^j$ wraps these annuli sufficiently large number of times around the $\breve{\phi}$-fibres in the direction opposite to the transverse orientation of $\breve{\mathcal{F}}_{1, s}$ as we pass from $\breve{\partial}_{1, s}^j$ to $\breve{\partial}_{2, s}^j$, $1\leqslant j\leqslant p+2$. Figure \ref{DT} shows one Dehn twist.

We perform Dehn twist operations on $\breve{\mathcal{F}}_{1,s}$, by two steps. At first, perform $\{D_s^j, 1\leqslant j\leqslant p\}$, then perform $D_s^{p+1}$ and $D_s^{p+2}$. Let  $\breve{\mathcal{F}}_{1, s}^{'}$ be the new surface bundle obtained by Dehn twist operations $\{D_s^j, 1\leqslant j\leqslant p\}$ on  $\breve{\mathcal{F}}_{1, s}$. By Remark \ref{intersection}, there are two cases for the foliation of $\breve{U}_{i,k,s}^{0,2l-1}$ by $\breve{U}_{i,k,s}^{0,2l-1}\cap\breve{\mathcal{F}}_{1, s}^{'}$, $s=1,2$.

\noindent Case A. $l=1$ and ${i-k}\leqslant\displaystyle{\frac{p-1}{2}}$, or $l=2$ and ${i-k}>\displaystyle{\frac{p-1}{2}}$. $1\leqslant i,k\leqslant p, i\neq k$

In this case $\overset{p}{\underset{j=1}{\cup}}N(\breve{V}_s^j)\cap \breve{L}_{i, k, s}^{2l-1}$ and $\overset{p}{\underset{j=1}{\cup}}N(\breve{V}_s^j)\cap \breve{L}_{i, k, s}^{2l+5}$ are both \textit{two} arcs. In particular $\breve{L}_{i,k,s}^{1}$ and $\breve{L}_{i,k,s}^7$ intersect $\breve{V}_s^i$ and $\breve{V}_s^k$, and $\breve{L}_{i,k,s}^{3}$ and $\breve{L}_{i,k,s}^9$ intersect $\breve{V}_{s}^k$ twice by Table \ref{orientL}. Then $\breve{L}_{i, k, s}^{2l-1}$ intersect $\overset{p}{\underset{j=1}{\cup}}N(\breve{V}_s^j)$ two times, and in different signs by Remark \ref{intersection}. If the tail of one arc of $\breve{L}_{i, k, s}^{2l-1}\cap \overset{p}{\underset{j=1}{\cup}}N(\breve{V}_s^j)$ lies in $\overset{p}{\underset{j=1}{\cup}}\breve{\partial}_{1, s}^j$, the tail of the other one lies in $\overset{p}{\underset{j=1}{\cup}}\breve{\partial}_{2, s}^j$. The same thing happens for $\breve{L}_{i, k, s}^{2l+5}$.

$\breve{U}_{i, k, s}^{0, 2l-1}\cap \breve{\mathcal{F}}_{1, s}^{'}$ is obtained from $\breve{U}_{i, k, s}^{0, 2l-1}\cap \breve{\mathcal{F}}_{1, s}$ by $D_s^j$'s ($1\leqslant j\leqslant p$, $s=1,2$) in a regular neighborhood of \textit{two} $\breve{\phi}$-circle fibres in $\breve{U}_{i, k, s}^{0, 2l-1}$ in the direction opposite to the transverse orientation of $\breve{\mathcal{F}}_{1, s}$ as we pass from $(\overset{p}{\underset{j=1}{\cup}}\breve{\partial}_{1, s}^j)\cap \breve{U}_{i, k, s}^{0, 2l-1}$ to $(\overset{p}{\underset{j=1}{\cup}}\breve{\partial}_{2, s}^j)\cap \breve{U}_{i, k, s}^{0, 2l-1}$, $1\leqslant l\leqslant 3$. Figure \ref{dehntwistU}-2 shows $\breve{U}_{i,k,s}^{0,3}\cap \breve{\mathcal{F}}_{1, s}^{'}$, and $\breve{U}_{i,k,s}^{0,1}\cap \breve{\mathcal{F}}_{1, s}^{'}$ is similar to $\breve{U}_{i,k,s}^{0,3}\cap \breve{\mathcal{F}}_{1, s}^{'}$.

Now we perform an isotope of $\breve{\mathcal{F}}_{1, s}^{'}$ in a small regular neighborhood of $\breve{U}_{i, k, s}^{0, 2l-1}$ in $\breve{Y}_{1, s}$, say $N(\breve{U}_{i, k, s}^{0, 2l-1})$, such that $\breve{U}_{i, k, s}^{0, 2l-1}\cap \breve{\mathcal{F}}_{1, s}^{'}$ changes back to $\breve{U}_{i, k, s}^{0, 2l-1}\cap \breve{\mathcal{F}}_{1, s}$ as shown in Figure \ref{dehntwistU}-1, and $\partial N(\breve{U}_{i, k, s}^{0, 2l-1})\cap \breve{\mathcal{F}}_{1, s}^{'}$ always keeps to be the same. The isotope is described as following. Push the whole arc $AB$ (like a finger move) along $\breve{\phi}$ in $\breve{U}_{i, k, s}^{0, 2l-1}$ in the direction opposite to the Dehn twist and the same times as the Dehn twist, meanwhile keep $\partial N(\breve{U}_{i, k, s}^{0, 2l-1}\cap \breve{\mathcal{F}}_{1, s}^{'})$ all the time (cf. Figure \ref{dehntwistU}-2). We still call the surface bundle $\breve{\mathcal{F}}_{1, s}^{'}$ after this isotope.

In case A, $\breve{U}_{i, k, s}^{0, 2l-1}\cap \breve{\mathcal{F}}_{1, s}^{'}$ is same as $\breve{U}_{i, k, s}^{0, 2l-1}\cap \breve{\mathcal{F}}_{1, s}$.

\noindent Case B. $l=3, 1\leqslant i, k\leqslant p, i\neq k$.

In this case, $\overset{p}{\underset{j=1}{\cup}}N(\breve{V}_s^j)\cap \breve{L}_{i, k, s}^{2l-1}$  and $\overset{p}{\underset{j=1}{\cup}}N(\breve{V}_s^j)\cap \breve{L}_{i, k, s}^{2l+5}$ are both \textit{one} arc. The foliation of $\breve{U}_{i, k, s}^{0, 5}$ determined by $\breve{U}_{i, k, s}^{0, 5}\cap \breve{\mathcal{F}}_{1, s}^{'}$ is obtained from the foliation $\breve{U}_{i, k, s}^{0, 5}\cap \breve{\mathcal{F}}_{1, s}$ by $D_s^j$'s $(1\leqslant j\leqslant p)$ in a regular neighborhood of one $\breve{\phi}$-circle fibre in $\breve{U}_{i, k, s}^{0, 5}$ in the direction opposite to the transverse orientation of $\breve{\mathcal{F}}_{1, s}$ as we pass from $(\overset{p}{\underset{j=1}{\cup}}\breve{\partial}_{1, s}^j)\cap \breve{U}_{i, k, s}^{0, 5}$ to$(\overset{p}{\underset{j=1}{\cup}}\breve{\partial}_{2, s}^j)\cap \breve{U}_{i, k, s}^{0, 5}$ (cf. Figure \ref{dehntwistU}-3).

Now we perform the Dehn twist operations $D_s^{p+1}$, $D_s^{p+2}$ on $\mathcal{F}'_{1,s}$.

Suppose $\breve{\mathcal{F}}_{1, s}^{''}$ is the new surface bundle of $\breve{Y}_{1, s}$ obtained by $D_s^{p+1}$ and $D_s^{p+2}$. By Remark \ref{intersection2}, $(N(\breve{V}_s^{p+1})\cup N(\breve{V}_s^{p+2}))\cap \breve{L}_{i, k, s}^{2l-1}$ and $(N(\breve{V}_s^{p+1})\cup N(\breve{V}_s^{p+2}))\cap \breve{L}_{i, k, s}^{2l+5}$ are both one arc, $l=1,2$. Then the situation of $\breve{U}_{i, k, s}^{0, 2l-1}\cap \breve{\mathcal{F}}_{1, s}^{''}$ is the same as case B in the first step (cf. Figure \ref{dehntwistU}-3), $l=1, 2$. 
The foliation of $\breve{U}_{i, k, s}^{0, 5}$ determined by $\breve{U}_{i, k, s}^{0, 5}\cap \breve{\mathcal{F}}_{1, s}^{''}$ is the same as $\breve{U}_{i, k, s}^{0, 5}\cap \breve{\mathcal{F}}_{1, s}^{'}$ because $(N(\breve{V}_s^{p+1})\cup N(\breve{V}_s^{p+2}))\cap(\breve{U}_{i,k,s}^{0,5})=\emptyset$.

Finally like \cite{abz}, we adjust $\breve{\mathcal{F}}_{1, s}^{''}$ by isotope which is the identity in a small regular neighborhood of $\partial \breve{Y}_{1, s}$ and outside a small regular neighborhood of \linebreak $\overset{3}{\underset{l=1}{\cup}}(\underset{1\leqslant i, k\leqslant p, i\neq k}{\cup}\breve{U}_{i, k, s}^{0, 2l-1})$ such that the interval foliation in each $\breve{U}_{i, k, s}^{0, 2l-1}$ becomes transverse to $\breve{L}_{i, k, s}^{2l-1}\cup \breve{L}_{i, k, s}^{2l+5}$, $1\leqslant l\leqslant 3$ (cf. Figure \ref{dehntwistU}-4). We denote the resulting surface bundle $\breve{\mathcal{F}}_{1, s}^{'''}$. Since $\breve{\mathcal{F}}_{1, s}^{'''}\cap \partial \breve{Y}_{1, s}$ is same as $\breve{\mathcal{F}}_{1, s}\cap \partial \breve{Y}_{1, s}$, the resulting surface bundle $\breve{\mathcal{F}}^{'''}$ in $\breve{M}$ is transverse to the link $\overset{2}{\underset{s=1}{\cup}}(\underset{1\leqslant i, k\leqslant p; i\neq k}{\cup}\breve{L}_{i, k, s})$ every where. In other words, the exterior of $\overset{2}{\underset{s=1}{\cup}}(\underset{1\leqslant i, k\leqslant p; i\neq k}{\cup}\breve{L}_{i, k, s})$ in $\breve{M}$ has a surface bundle structure, i.e., $\overset{2}{\underset{s=1}{\cup}}(\underset{1\leqslant i, k\leqslant p; i\neq k}{\cup}\breve{L}_{i, k, s})$ is a fibred link, i.e. $K$ is a virtually fibred Montesinos knot.  Now we finish the proof of Theorem \ref{main} when $K$ is a classic Montesinos knot, and $n=3$.

\subsection{$\mathbf{K}$ is a link of two components.}

 The proof in this case is mostly same as in section 2.1, except for the following changes.

In this case, $f |: \widetilde{K}\rightarrow K^*$ is a trivial 2-fold cover. Then $L=\Psi^{-1}(\widetilde{K})$ is a geodesic link with $2p^2$ components. Let $L=\{L_{i,j}^r: 1\leqslant i,j \leqslant p, r=1,2\}$. $\hat{f}(L_{i,j}^r)=L_{i,j}^*$, $r=1,2$.
 Now we take $M$ to be the complement  of $\{L_{j, j}^1: 1\leqslant j \leqslant p\}$ in $Y$. As in the proof of Proposition \ref{semibundle}, we have $M=Y_1\cup M_2^1\cup M_2^2\cup\cdots\cup M_2^p$, where $M_2^j=\hat{f}^{-1}(L_{j,j}^*\times[-\epsilon, \epsilon])-\overset{\circ}{N}(L_{j,j}^1)$, $1\leqslant j\leqslant p$. Note that $L_{j, j}^2$ is a fibre of the Seifert fibre structure on $M_2^j, 1\leqslant j\leqslant p$.

Since $K$ has two components, and $p$ is odd, we have $q_1+q_2+q_3$ is even.
\begin{equation*}
e=p^2e(W_K)=-p(q_1+q_2+q_3),\ \ \ \ \ \ \tilde{e}=\frac{e}{p}=-(q_1+q_2+q_3)
\end{equation*}
so $|\tilde{e}|$ is even and non-zero. We use the definitions of $a, b, c, d$ in the proof of Proposition \ref{semibundle} , and in this case $c=1$. We can take $a=c=d=1$ and $b=0$. i.e.
\begin{equation*}
\overline{\alpha}^j=\alpha^j, \ \ \ \ \ \ \overline{\phi}^j=\alpha^j+\phi^j,\ \ \ \ \ \ 1\leqslant j\leqslant p.
\end{equation*}

In the following discussion, we will always assume that $1\leqslant j\leqslant p$ except for special indication.

By a similar analysis in Proposition \ref{semibundle}, we can see that there are orientable horizontal surfaces $H_1$ in $Y_1$, and $H_2^j$ in $M_2^j$ which piece together to form the semifibre $H$ where the projection of $H_1$ to $F_1$ has degree $|\lambda|=1$ and that of $H_2^j$ to the base of $M_2^j$ has degree $|\overline{\lambda}|=\displaystyle|\frac{\tilde{e}}{2}|$.

Further, if we suppose that the slope of $H_1$ on $T_i^j$ is given by $u_i^j\alpha_{1, i}^j+t_i^j\phi_{1, i}^j\ (i=1, 2)$, and that of $H_2^j$ on $T_i^j$ is $\overline{u}_i^j\alpha_{2, i}^j+\overline{t}_i^j\phi_{2, i}^j\ (i=1, 2, 3)$. Then we get the following result from (\ref{slope12}) and (\ref{slope3})

\begin{equation*}
\frac{t_1^j}{u_1^j}=\frac{\tilde{e}-2}{2},\  \frac{t_2^j}{u_2^j}=\frac{2-\tilde{e}}{2}; \ \ \ \ \ \ \ \ \frac{\overline{t}_1^j}{\overline{u}_1^j}=\frac{2-\tilde{e}}{\tilde{e}}, \ \frac{\overline{t}_2^j}{\overline{u}_2^j}=\frac{2+\tilde{e}}{\tilde{e}}, \ \frac{\overline{t}_3^j}{\overline{u}_3^j}=-\frac{4}{\tilde{e}}.
\end{equation*}

Since $\tilde{e}$ is even, the values of the coefficients $u_i^j, t_i^j, \overline{u}_i^j, \overline{t}_i^j$ are given as following
\begin{align*}
t_1^j=\frac{\tilde{e}}{2}-1,\ \  u_1^j=1; \ \ \ \ \ \   t_2^j=1-\frac{\tilde{e}}{2},\ \ u_2^j=1;\\
\overline{t}_1^j=1-\frac{\tilde{e}}{2},\ \ \overline{u}_1^j=\frac{\tilde{e}}{2}; \ \ \ \ \  \overline{t}_2^j=1+\frac{\tilde{e}}{2}, \ \ \overline{u}_2^j=\frac{\tilde{e}}{2}; \\
(\overline{t}_3^j , \overline{u}_3^j)=
\begin{cases}
 (-2, \displaystyle{\frac{\tilde{e}}{2}})& \text{if $\tilde{e}=4k+2$},\\[2ex]
 (-1, \displaystyle{\frac{\tilde{e}}{4}})& \text{if $\tilde{e}=4k$},
\end{cases}
k\in\mathbb{Z}.
\end{align*}

As in 2.1 $H_2^j$ is a surface which interpolates between the slope $\displaystyle{\frac{2-\tilde{e}}{\tilde{e}}}$ on $T_1^j$ and $\displaystyle{\frac{2+\tilde{e}}{-\tilde{e}}}$ on $T_2^j$. $\partial H_1\cap T_{1, r}^j$ has $\displaystyle{|\frac{\lambda}{u_r^j}|=|\frac{1}{1}|=1}$ component $r=1, 2$, and $\partial H_2^j\cap T_{2, r}^j$ has $\displaystyle{|\frac{\overline{\lambda}}{\overline{u}_r^j}|=|\frac{\tilde{e}/2}{\tilde{e}/2}|=1}$ component $r=1, 2, 3$. Still $\epsilon_1=-\epsilon_2$ implies that $M$ is a surface semi-bundle. Denote the associated fibring in $M_2^j$ by $\mathcal{F}_2^j$.  $\mathcal{F}_2^j$ is transverse to all the new fibres, and in particular to $L_{j, j}^2$.

Now $L_{j, j}^1$ only intersects the original fibre once, so $U_{i, j}^{0, 2l}$ is a once-punctured annulus, and there is only one singular point on $U_{i, j}^{0, 2l}, 1\leqslant i, j\leqslant p, i\neq j, l=1,2,3$.

\begin{figure}
\begin{center}
\includegraphics[width=5in]{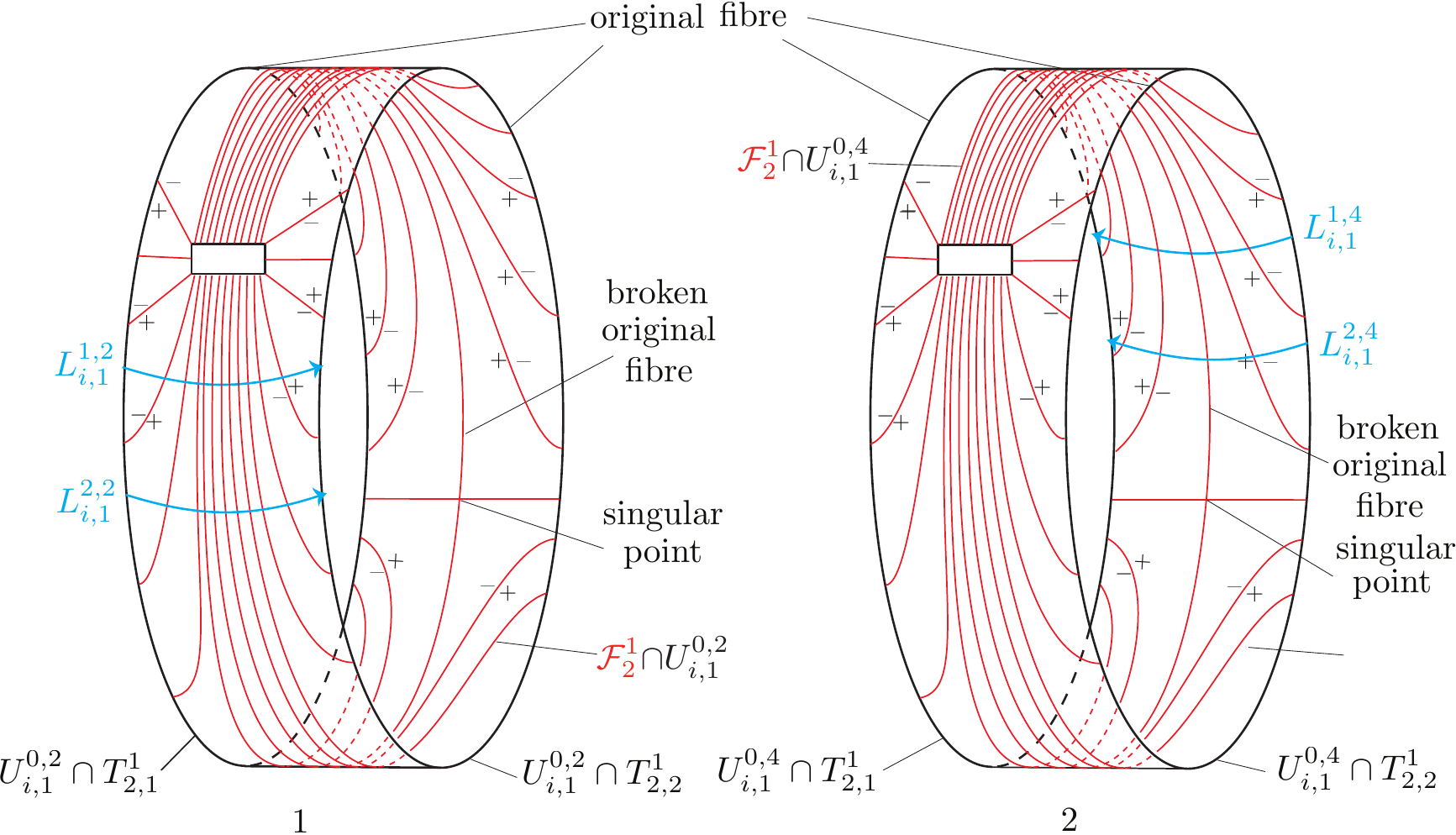}
\end{center}
\caption{\label{foliation2}$U_{i, 1}^{0, 2}$ and $U_{i, 1}^{0, 4}$.}
\end{figure}

The intersection of $U_{i, 1}^{0, 2}$ with $\mathcal{F}_2^1$ is depicted in Figure \ref{foliation2}-1,  $U_{i, 1}^{0, 4}\cap \mathcal{F}_2^1$ is depicted in Figure \ref{foliation2}-2, and $U_{1, i}^{0, 6}\cap \mathcal{F}_2^1$ is similar to $U_{i, 1}^{0, 4}\cap \mathcal{F}_2^1, 2\leqslant i\leqslant p$. The previous argument can be proceed without significant change to produce the desired surface bundle structure on $\breve{M}_{2, s}^j$ which is transverse to $\breve{L}'=\underset{s=1}{\overset{2}{\cup}}(\underset{r=1}{\overset{2}{\cup}}\underset{1\leqslant i, j\leqslant p; i\neq j}{\cup}\breve{L}_{i, j, s}^r\cup(\overset{p}{\underset{j=1}{\cup}}\breve{L}_{j, j, s}^2))$. By similar Dehn twist operation on $\breve{\mathcal{F}}_{1, s}$, we can get a surface bundle in $\breve{Y}_{1, s}$, say $\breve{\mathcal{F}}_{1, s}^{'}$ such that $\breve{\mathcal{F}}_{1, s}^{'}$ is transverse to $\breve{L}', s=1, 2$. The proof of Theorem \ref{main} is finished in Case (1).

\section{Proof of Theorem \ref{main} in Case (2).}

In this section, we consider classic Montesinos link $K( {q_1}/{p}, {q_2}/{p}, \cdots, {q_n}/{p})$, where $n>3$, $p$ odd and $p\geq 3$.  We do not consider the case $e(W_K)=0$, otherwise $K$ is virtually fibred by Theorem 1.2 \cite{abz}.

We still consider $K$ is a knot first, and prove Theorem \ref{main} by the same four steps as in Sec 2.1.

The Seifert structure is $f: W_K\rightarrow S^2(\underbrace{p, p, \cdots, p}_{n})$, and $\Gamma_1=\pi_1(S^2(\underbrace{p, p, \cdots, p}_{n}))$ has a presentation
\begin{equation*}
\Gamma_1=<x_1, x_2, \cdots, x_n\ |\ x_1^p= \cdots =x_n^p=x_1x_2 \cdots x_n=1>.
\end{equation*}
Define $h_1$ similarly as in Sec. 2.1.
\begin{equation*}
h_1: \Gamma_1\rightarrow \mathbb{Z}/p
\end{equation*}
where $h_1(x_1)=\overline{1}, h_1(x_2)=-\overline{1}$, and $h_1(x_i)=\overline{0}$, for $3\leqslant i\leqslant n$. Let $\psi_1: F'\rightarrow S^2(\underbrace{p, p, \cdots, p}_{n})$ be the $p$-fold cyclic orbifold cover of $S^2(\underbrace{p, p, \cdots, p}_{n})$ corresponding to $h_1$. Then $F'=S^2(\underbrace{p, p, \cdots, p}_{(n-2)p})$ and
\begin{equation*}
\Gamma_2=\pi_1(F')=<y_1, y_2, \cdots, y_{(n-2)p}\ |\ y_1^p=y_2^p= \cdots =y_{(n-2)p}^p=y_1y_2 \cdots y_{(n-2)p}=1>.
\end{equation*}
Define
\begin{equation*}
h_2: \Gamma_2\rightarrow \mathbb{Z}/p
\end{equation*}
where $h_2(y_i)=\overline{1}$ for $1\leqslant i\leqslant (n-2)p$. Suppose that $\psi_2: F\rightarrow F'$ is the $p$-fold cyclic orbifold cover of $F'$ corresponding to $h_2$, and $\psi=\psi_1\circ \psi_2$. Now let $\Psi: Y\rightarrow W_K$ be the associated $p^2$-fold cover. $L=\Psi^{-1}(\widetilde{K})$ is still a geodesic link in $Y$ with exactly $p^2$ components. We use the same notations as in section 2. $L=\{L_{i, j}, 1\leqslant i,j\leqslant p\}$, $\hat{f}(L_{i, j})=L_{i, j}^*$, $1\leqslant i, j\leqslant p$, and $\psi^{-1}(c_i)=\hat{c}_{i, j}$, $1\leqslant i\leqslant n, 1\leqslant j\leqslant p$.

Let $M=Y\setminus (\overset{p}{\underset{j=1}{\cup}}\overset{\circ}{N}(L_{j, j}))$. By the same process from Step 2 to Step 4 as described in section 2.1, we can find a double cover of $M$, say $\breve{M}$ and a surface bundle $\breve{\mathcal{F}}^{''}$ in $\breve{M}$ such that it is transverse to the lift of $L_{i, j}$, $1\leqslant i, j\leqslant p, i\neq j$.

The only thing different here is the construction of $\Gamma$. In Step 4, we perform Dehn twist operations on $\breve{\mathcal{F}}_{1, s}$ along a set of $\breve{\phi}$-vertical tori $\Gamma$ such that the new surface fibre is transverse to $\breve{L}_{i, k, s}^{2l-1}$ everywhere, $1\leqslant l\leqslant 2n$, $s=1,2$, $1\leqslant i,k\leqslant p$, $i\neq k$. $\Gamma=\{\breve{V}_s^1,\breve{V}_s^2 : s=1,2\}$ contains four tori not $2p+4$ as in Case (1). They are the lifts of two $\phi$-vertical tori, $V^r=\hat{f}^{-1}(l_r)$, $r=1, 2$. $l_1$ and $l_2$ are two simple closed curves on $F$, which intersect all $(L_{i, k}^{2l-1})^*$'s once or twice when $n\geq 5$ in the same signs, $1\leqslant i, k\leqslant p, i\neq k$, and $1\leqslant l\leqslant n$.

Note that we always suppose that $1\leqslant i,j,k \leqslant p$, $i\neq k$

Recall that $F_2^j=L_{j,j}^*\times[-\epsilon, \epsilon]$, and $F_1= F-\overset{p}{\underset{j=1}{\cup}}{\overset{\circ}{F_2^j}}$. As before, $L_{i,k}^*$ is separated into $2n$ parts by $\overset{p}{\underset{j=1}{\cup}}F_2^j$ and $F_1$. Table \ref{L*} is still true in this case. In addition, if $4\leqslant l\leqslant n$,

\begin{center}
the tail of $(L_{i,k}^{2l-1})^*$ in
$\begin{cases}
 \beta_2^i     & \displaystyle{i-k}\leqslant \frac{p-1}{2}, \\[2ex]
  \beta_1^i    & \displaystyle{i-k}>\frac{p-1}{2},
\end{cases}
$

the head of $(L_{i,k}^{2l-1})^*$ in
$\begin{cases}
 \beta_1^i     & \displaystyle{i-k}\leqslant \frac{p-1}{2}, \\[2ex]
  \beta_2^i    & \displaystyle{i-k}>\frac{p-1}{2}.
\end{cases}$
\end{center}

As in Step 3 in Section 2, we can \textit{reorient} $\{L_{i,k}, 1\leqslant i,k \leqslant p, i\neq k\}$ such that

\begin{center}
the tail of
$\begin{cases}
  \text{$L_{i,k}^2$  and  $L_{i,k}^{2+2n}$ lie in $T_1^k$}, \\
  \text{$L_{i,k}^4$  and  $L_{i,k}^{4+2n}$ lie in $T_2^k$}, \\
  \text{$L_{i,k}^{2l}$  and  $L_{i,k}^{2l+2n}$ lie in $T_2^i$, $3\leqslant l\leqslant n$}.
\end{cases}$
\end{center}

Table \ref{orientL} is also true in this case and the tail of $L_{i,k}^{2l-1+r}$ is in $T_1^i$, and the head of $L_{i,k}^{2l-1+r}$ is in $T_2^i$, $4\leqslant l\leqslant n$, $r=0,2n$.

The construction of $l_1$ and $l_2$ is similar as before.

$a_r^j$ is still the intersection of $\beta_r^j$ and a small disk centered at $\hat{c}_{1,j}$ with radius bigger than $\epsilon$, $1\leqslant j\leqslant p$, $r=1,2$. Push $a_r^j$ into $F_1$ a little bit. $b_r^j$ is different.  We take $b_r^j$ in $F_1$ which is parallel to $\beta_r^j$ and between $\hat{c}_{2,j}$ and $\hat{c}_{n,j}$, $1\leqslant j\leqslant p$, $r=1,2$.
Then  $\overset{p}{\underset{j=1}{\cup}}(a_1^j\cup a_2^j)$ intersects $(L_{i,k}^{2l-1})^*$ once, where $l=1,2$, and $\overset{p}{\underset{j=1}{\cup}}(b_1^j\cup b_2^j)$,  intersects $(L_{i,k}^{5})^*$ and $(L_{i,k}^{2n-1})^*$ both once and intersect $(L_{i,k}^{2l-1})^*$ twice where $3<l<n$, $1\leqslant i,k\leqslant p$, $i\neq k$.

\begin{figure}
\begin{center}
\includegraphics[width=5.9in]{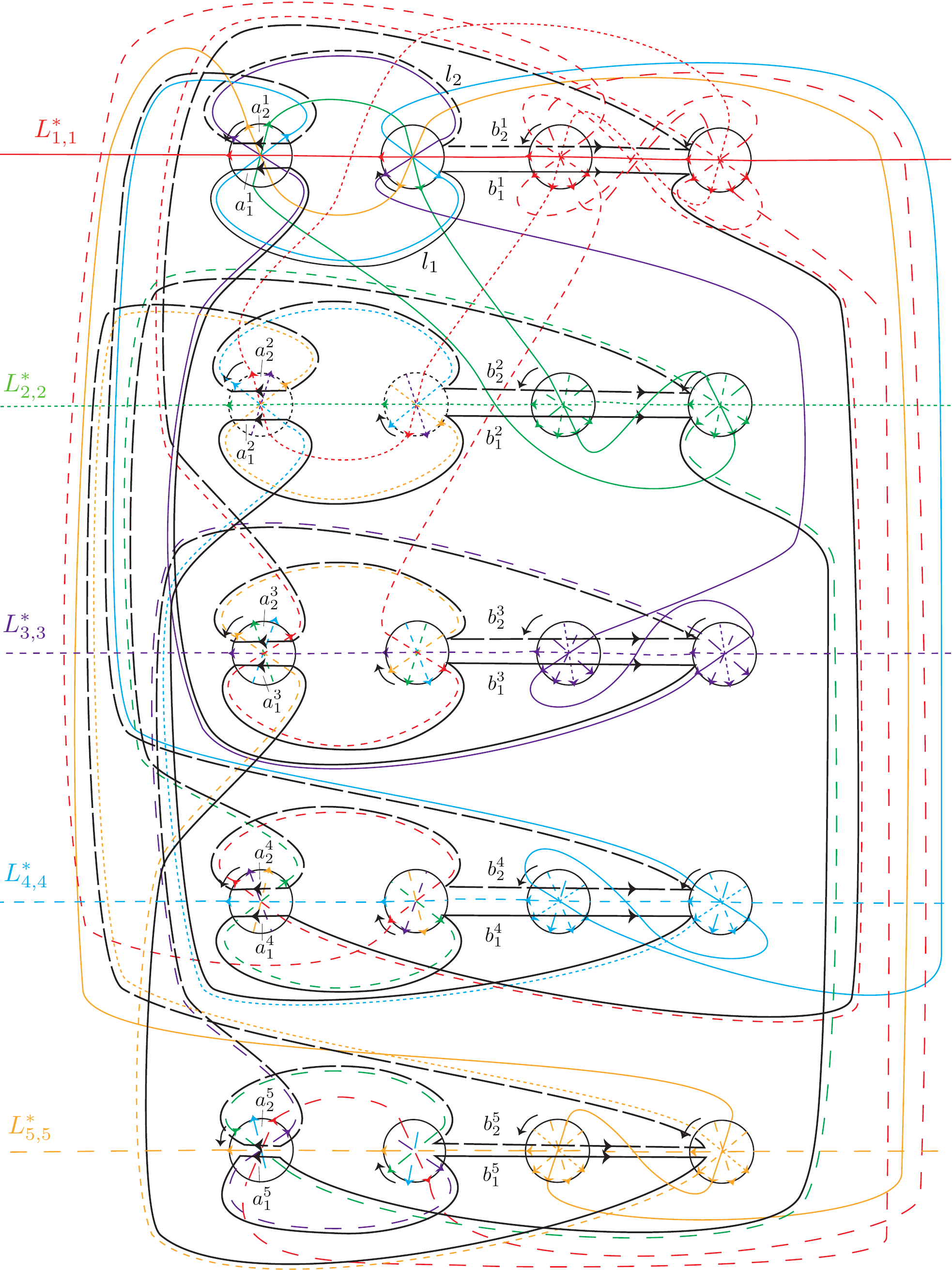}
\end{center}
\caption{\label{n=4}$l_1$ and $l_2$ when $n=4$.}
\end{figure}

\begin{figure}
\begin{center}
\includegraphics[width=5.9in]{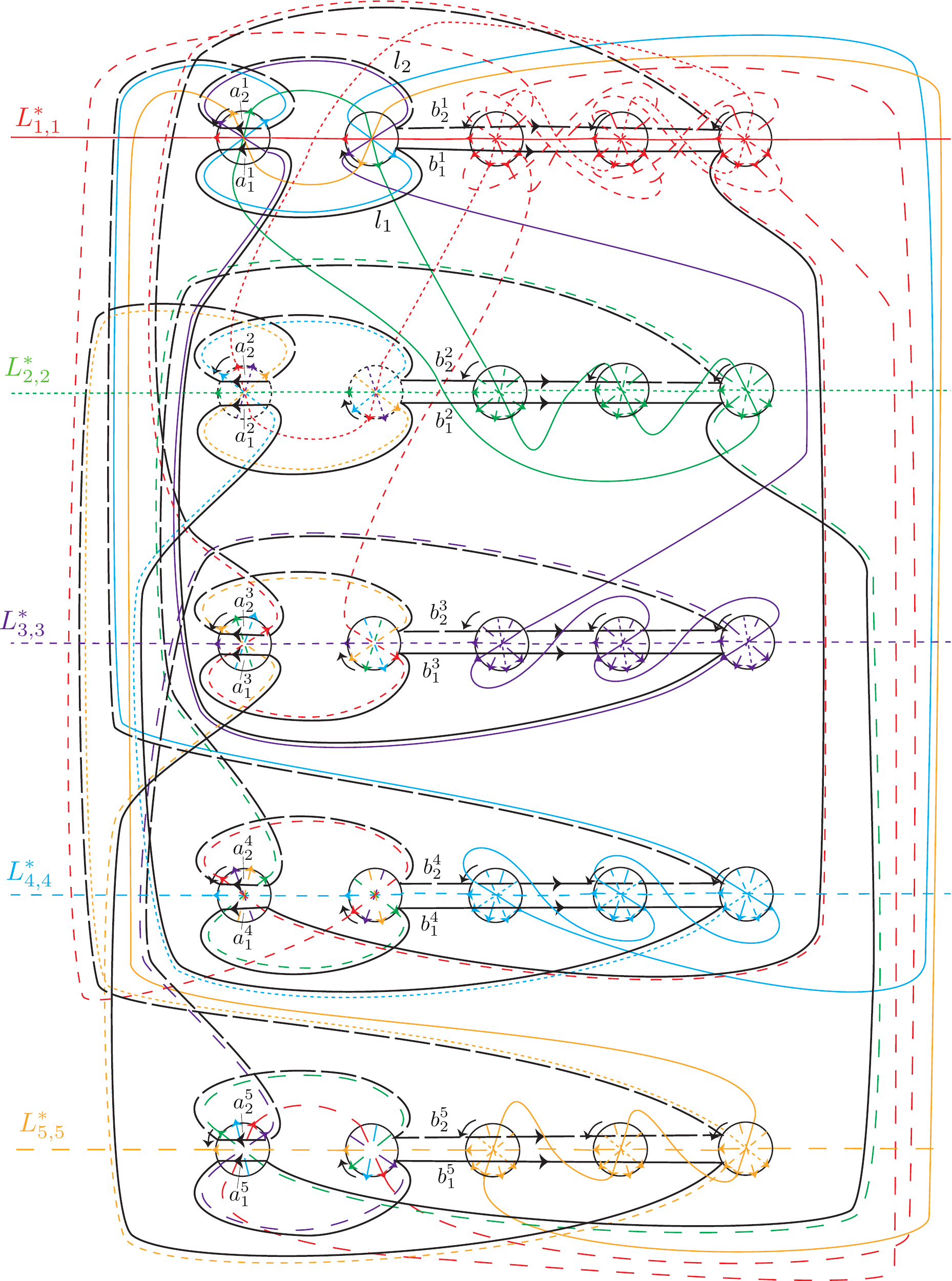}
\end{center}
\caption{\label{n=5}$l_1$ and $l_2$ when $n=5$.}
\end{figure}

We still can connect $a_r^j$'s and $b_r^j$'s by the same arcs as in (\ref{connect}), $1\leqslant j\leqslant p$, $r=1,2$. See Figure \ref{n=4} in the case $n=4$, and Figure \ref{n=5} in the case $n=5$. When $n>5$, $l_1$ and $l_2$ are similar to the ones in the case $n=5$.


 $\breve{V}_s^1=p_2^{-1}(V^1)\cap \breve{Y}_{1, s}$, $\breve{V}_s^2=p_2^{-1}(V^2)\cap \breve{Y}_{1, s}$, $s=1, 2$. Recall that $p_2: \breve{M}\rightarrow M$ is the double cover of $M$.
We have Remark \ref{intersection3} by the construction of $l_1$ and $l_2$.

\begin{remark}\label{intersection3}
$\breve{L}_{i,k,s}^{2l-1+r}$'s intersect $\{\breve{V}_s^1, \breve{V}_s^2: s=1,2\}$ negatively. In addition, $\breve{L}_{i,k,s}^{2l-1+r}$ intersects $\{\breve{V}_s^1, \breve{V}_s^2: s=1,2\}$ exactly once when $l=1,2,3,n$, and twice when $3<l<n$, where $1\leqslant i,k \leqslant p$, $i\neq k$, $r=0,2n$.
\end{remark}

By Remark \ref{intersection3} we can define $\partial(N(\breve{V}_s^r))=\breve{\partial}_{1, s}^r\cup \breve{\partial}_{2, s}^r$ such that the tail of $\breve{L}_{i, k, s}^{2l-1}$ and $\breve{L}_{i, k, s}^{2l-1+2n}$ on $\breve{\partial}_{1, s}^r$, $r=1, 2$, $s=1, 2$, $1\leqslant i, k\leqslant p$, $i\neq k$, $1\leqslant l\leqslant n$.

\begin{figure}
\begin{center}
\includegraphics[width=4.5in]{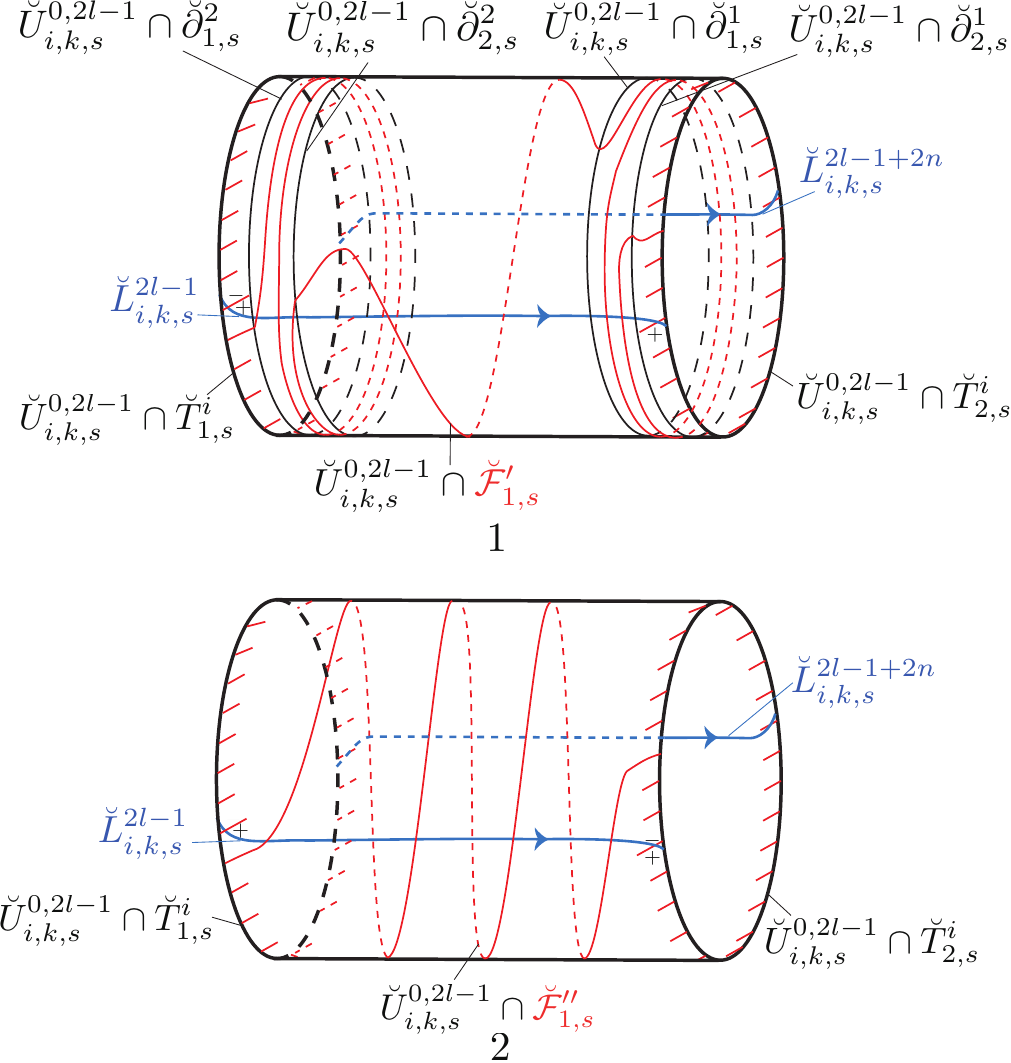}
\end{center}
\caption{\label{foliationU2} Dehn twist of $\breve{\mathcal{F}}_{1, s}$.}
\end{figure}

Now we perform a Dehn twist operation $D_s^r$ which wraps these annuli $m$ times around the $\breve{\phi}$-fibres in the direction opposite to the transverse orientation of $\breve{\mathcal{F}}_{1, s}$ as we pass from $\breve{\partial}_{1, s}^r$ to $\breve{\partial}_{2, s}^r$, $r=1,2$, $s=1,2$.
By Remark \ref{intersection3} $\breve{L}_{i, k, s}^{2l-1}$ intersects $(\breve{V}_s^1\cup \breve{V}_s^2)$ once, $\breve{U}_{i, k, s}^{0, 2l-1}\cap \breve{\mathcal{F}}_{1, s}^{'}$ is similar as in Figure \ref{dehntwistU}-3, when $l=1,2,3,n$. When $3<l<n$, $\breve{L}_{i, k, s}^{2l-1}$ intersects $(\breve{V}_s^1\cup \breve{V}_s^2)$ twice and in the same sign. $\breve{U}_{i, k, s}^{0, 2l-1}\cap \breve{\mathcal{F}}_{1, s}^{'}$ is shown in Figure \ref{foliationU2}-1.
In this case, $\breve{L}_{i, k, s}^{2l-1}$ goes from $\breve{T}_{1, s}^i$ to $\breve{T}_{2, s}^i$. Now $\breve{U}_{i, k, s}^{0, 2l-1}\cap \breve{\mathcal{F}}_{1, s}^{'}$ is wrapped $2m$ times around the $\breve{\phi}$-fibres in the direction opposite to the transverse orientation of $\breve{\mathcal{F}}_{1, s}$ as we pass from $\breve{U}_{i, k, s}^{0, 2l-1}\cap \breve{T}_{1, s}^{i}$ to $\breve{U}_{i, k, s}^{0, 2l-1}\cap \breve{T}_{2, s}^{i}$, $1\leqslant i,k\leqslant p$, $i\neq k$, $3<l<n$. We adjust $\breve{\mathcal{F}}_{1, s}^{'}$ by isotope, and denote the resulting surface bundle $\breve{\mathcal{F}}^{''}$ in $\breve{M}$ such that it is transverse to $\breve{L}_{i, k, s}^{2l-1}$ as shown in Figure \ref{foliationU2}-2.

Now we finish the proof of Theorem \ref{main} in Case (2) when $K$ is a knot. If $K$ has two components, Theorem \ref{main} can be proved by combining Sec. 2.2 with the argument in this section.

\end{document}